\documentclass[letterpaper, 11pt]{amsart}

\usepackage{amsmath,amsthm,amsfonts,amssymb,amscd}
\usepackage{bbm}
\usepackage{bm}
\usepackage{tikz}
\usepackage{tikz-cd}
\usepackage{appendix}
\usepackage{BOONDOX-calo}
\usepackage{standalone}
\usepackage[hidelinks]{hyperref}
\usetikzlibrary{arrows,chains,matrix,positioning,scopes}
\usepackage[letterpaper, left=3cm,right=3cm, top=3cm, bottom=3cm]{geometry}
\usepackage{adjustbox}
\usepackage{enumitem}

\newcommand{\CC}{{\mathbb C}}

\newcommand{\ZZ}{{\mathbb Z}}
\newcommand{\QQ}{{\mathbb Q}}
\newcommand{\Ext}{{\textup{Ext}}}
\newcommand{\Extpan}{{\textup{Extpan}}}
\newcommand{\Hom}{{\textup{Hom}}}
\newcommand{\End}{{\textup{End}}}
\newcommand{\GL}{{\mathrm{GL}}}
\newcommand{\Gr}{{\mathrm{Gr}}}
\newcommand{\Id}{{\mathrm{Id}}}
\newcommand{\ab}{{\mathrm{ab}}}
\newcommand{\inEnd}{{\underline{\textup{End}}}}
\newcommand{\fu}{{\mathfrak{u}}}
\newcommand{\fv}{{\mathfrak{v}}}
\newcommand{\fg}{{\mathfrak{g}}}
\newcommand{\fR}{{\mathfrak{R}}}
\newcommand{\inHom}{{\underline{\textup{Hom}}}}
\newcommand\dual{\raise0.9ex\hbox{$\scriptscriptstyle\vee$}}
\newcommand{\sN}{{\mathcal{N}}}
\newcommand{\sL}{{\mathcal{L}}}
\newcommand{\sM}{{\mathcal{M}}}
\newcommand{\sE}{{\mathcal{E}}}

\newcommand{\bT}{{\mathbf{T}}}
\newcommand{\FF}{{\mathbb{F}}}
\newcommand{\bS}{{\mathbf{S}}}

\theoremstyle{plain}
\newtheorem{thm}{Theorem} 
\newtheorem{prop}[thm]{Proposition}
\newtheorem{lemma}[thm]{Lemma}  
\newtheorem{cor}[thm]{Corollary}
\numberwithin{thm}{subsection}

\newtheorem{manualtheoreminner}{Theorem}
\newenvironment{thm'}[1]{%
  \renewcommand\themanualtheoreminner{#1}%
  \manualtheoreminner
}{\endmanualtheoreminner}

\theoremstyle{definition}

\newtheorem{notation}[thm]{Notation}

\theoremstyle{remark}
\newtheorem{rem}[thm]{Remark}

\tikzset{>=stealth}

\makeatletter
\def\@seccntformat#1{%
  \protect\textup{\protect\@secnumfont
    \ifnum\pdfstrcmp{subsection}{#1}=0 \bfseries\fi
    \csname the#1\endcsname
    \protect\@secnumpunct
  }%
}  
\makeatother

\makeatletter
\@namedef{subjclassname@2020}{%
  $2020$ Mathematics Subject Classification}
\makeatother

\begin{document}

\title[Tannakian fundamental groups of blended extensions]{Tannakian fundamental groups of blended extensions}
\author{Payman Eskandari}
\address{Department of Mathematics and Statistics, University of Winnipeg, Winnipeg MB, Canada }
\email{p.eskandari@uwinnipeg.ca}
\begin{abstract}
Let $A_1, A_2,A_3$ be objects in a neutral tannakian category over a field of characteristic zero. Let $L$ be an extension of $A_2$ by $A_1$, and $N$ an extension of $A_3$ by $A_2$. Let $M$ be a blended extension ({\it extension panach\'{e}e}) of $N$ by $L$. We study the subgroup of the tannakian fundamental group of $M$ that acts trivially on $A_1, A_2,A_3$. We also give an application to the unipotent part of the motivic version of the Hodge conjecture (i.e., the equality of the unipotent radicals of the motivic Galois and Mumford-Tate groups) for Deligne 1-motives.
\end{abstract}
\maketitle

\section{Introduction}\label{sec:introduction}
The notion of a blended extension, invented by Grothendieck in \cite{Gr68}, provides a natural framework to study 3-step filtrations. By definition, given a fixed extension $L$ of $A_2$ by $A_1$ and a fixed extension $N$ of $A_3$ by $A_2$ in an abelian category, a blended extension of $N$ by $L$ is a diagram of the form 
\begin{equation}\label{eq1}
\begin{tikzcd}
   & & 0 \arrow{d} & 0 \arrow{d} &\\
   0 \arrow[r] & A_1 \ar[equal]{d} \arrow[r, ] & L  \arrow[d] \arrow[r, ] &  A_2 \arrow{d} \arrow[r] & 0 \\
   0 \arrow[r] & A_1 \arrow[r] & M \arrow[d] \arrow[r] &  N \arrow{d}  \arrow[r] & 0 \\
   & & A_3 \arrow{d} \ar[equal]{r} & A_3 \arrow{d} & \\
   & & 0 & 0 &   
\end{tikzcd}
\end{equation}
in the category with exact rows and columns. Here, with abuse of notation an extension and its middle object are denoted by the same letter, and the top row and right column are our two given fixed extensions $L$ and $N$.

In this paper we shall consider blended extensions $M$ as above in the setting of a neutral tannakian category $\mathbf{T}$ over a field of characteristic zero. Assume that $A_1$, $A_2$ and $A_3$ are semisimple objects. We are interested in the unipotent radical of the (tannakian) fundamental group of $M$. More explicitly, let $\fu(M)$ be the Lie algebra of the unipotent radical of the fundamental group of $M$; it is a canonical subobject of $\inEnd(M):= \inHom(M,M)$ (where $\inHom$ is the internal Hom), whose image under any fiber functor $\omega$ is the Lie algebra of the unipotent radical of the fundamental group of $M$ with respect to $\omega$. Our subject of study is $\fu(M)$, and we would like to describe it, ideally as explicitly and computably as possible, in terms of the extensions that appear in \eqref{eq1}.

The analogue of this problem for extensions, i.e. the determination of the unipotent radical of the fundamental group of an extension $L$ of a semisimple object $A_2$ by a semisimple object $A_1$, has been studied by Bertrand \cite{Ber01} in the special setting of differential equations and by Hardouin (\cite{Har05}, \cite{Har06}) and the author and Murty \cite{EM1} in general.

The problem of determination of $\fu(M)$ for blended extensions $M$ with semisimple  $A_1$, $A_2$ and $A_3$ arises naturally in at least two important settings. The first is in connection to motives with 3 weights (i.e. where the associated graded with respect to the weight filtration has 3 graded components) and their realizations (e.g. Hodge and $\ell$-adic realizations). The weight filtration on a motive with 3 weights gives rise to a blended extension as above in the motivic or realization category. Depending on whether one considers the blended extension in a tannakian category of motives or a tannakian category of realizations (e.g. rational mixed Hodge structures for the Hodge realization), $\fu(M)$ is the Lie algebra of the unipotent radical of the motivic Galois group of $M$ or the corresponding group\footnote{This assumes that the realization in question is such that the realizations of pure motives are semisimple. For the Hodge realization this is known, but it is conjectural for the $\ell$-adic and de Rham - Betti realizations.} for the choice of realization (e.g. the Mumford-Tate group in the case of Hodge realization).

The second natural setting where this problem arises is in connection to the differential Galois groups of products of three completely reducible differential operators. This second setting of the problem has been studied by Bertrand \cite{Ber01} and Hardouin \cite{Har05} (see also the references therein). In particular, according to \cite{Har05} one has a complete description of $\fu(M)$ in this setting.

Focusing on the setting of motives, arguably the most accessible nontrivial class of motives with three weights are those coming from Deligne's theory of 1-motives over a field \cite{De74}. In this case, the problem of determination of $\fu(M)$ is well-understood. Bertolin and her collaborators (in a series of papers, starting with \cite{Be02} and \cite{Be03}, and most recently in \cite{BP1}; see also the references therein) and Jossen \cite{Jo14} have independently studied $\fu(M)$ in this case, and give explicit descriptions of $\fu(M)$ both for the motivic Galois and Mumford-Tate groups. Jossen also does this for the $\ell$-adic realization, proving the Mumford-Tate conjecture for 1-motives at the unipotent level.

This paper came out of the author's attempt to better understand aspects of the works of Bertolin and Jossen, and curiosity to find out how much of their results could be obtained more abstractly, without using the explicit geometric situation particular to 1-motives. In fact, analogues to some of their results can be seen in works of Bertrand and Hardouin on differential equations, hinting that some of the results should indeed hold more generally. This connection was surely already observed by (at least) Bertrand, as it is clear from his paper \cite{Ber13} on self-dual blended extensions in an arbitrary tannakian category.

For convenience and to help the reader navigate through the paper with more ease, we have summarized below most of the results of the paper on the general description of $\fu(M)$. The reader familiar with the works cited earlier in the contexts of differential equations and 1-motives will be able to see the connections between the statements below and some of the results in those contexts. Note that the notations $\Hom$ and $\Ext^1$ (with no subscript) refer to the Hom and Ext groups for our category $\bT$, and $\inHom$ is the internal Hom. 

\begin{thm'}{A}\label{thm: intro}
Let $M$ be the blended extension \eqref{eq1} in a neutral tannakian category $\bT$ over a field of characteristic zero. Suppose $A_1$, $A_2$, $A_3$ are semisimple, and that
\begin{equation}\label{eq35}
\Hom(\inHom(A_2,A_1)\oplus \inHom(A_3,A_2), \inHom(A_3,A_1)) = 0.
\end{equation}
Then we have the following:
\begin{itemize}[wide]
\item[(a)] The object $\fu(M)$ fits in a canonical exact sequence
\begin{equation}\label{eq36}
\begin{tikzcd}
0 \arrow[r] & \fu_{-2}(M) \arrow[r] & \fu(M) \arrow[r] & \fu_{-1}(M) \arrow[r] & 0
\end{tikzcd}
\end{equation}
where $\fu_{-2}(M):= \fu(M)\cap \inHom(A_3,A_1)$ and $\fu_{-1}(M)$ is a canonical subobject of 
\[
V:=\inHom(A_2,A_1)\oplus \inHom(A_3,A_2).
\]
Moreover, the subobjects $\fu_{-2}(M)$ and $\fu_{-1}(M)$ of $\inHom(A_3,A_1)$ and $V$ completely determine the subobject $\fu(M)$ of $\inEnd(M)$. {\em (See \S \ref{sec: construction of pi} and \S \ref{sec: def of u-1 and u-2} for the construction of the sequence. See Proposition \ref{prop: when u is determined by its two subquotients}, Proposition \ref{prop: decomposition of u as superdiagonals plus top right entries} including its proof, and \S \ref{sec: conditions C1 and C2} for the determination of $\fu(M)$ from $\fu_{-2}(M)$ and $\fu_{-1}(M)$.)}
\item[(b)] (Determination of $\fu_{-1}(M)$) Let $\sL$ and $\sN$ be the elements of $\Ext^1(\mathbbm{1}, \inHom(A_2,A_1))$ and $\Ext^1(\mathbbm{1}, \inHom(A_3,A_2))$ corresponding to $L$ and $N$. Consider the element
\[(\sL,\sN)\in \Ext^1(\mathbbm{1},\inHom(A_2,A_1))\oplus \Ext^1(\mathbbm{1},\inHom(A_3,A_2)) \cong \Ext^1(\mathbbm{1},V).\]
Then $\fu_{-1}(M)$ is the smallest subobject of $V$ with the property that the pushforward of $(\sL,\sN)$ to an extension of $ \mathbbm{1}$ by $V/\fu_{-1}(M)$ splits. In particular, the subobject $\fu_{-1}(M)$ of $V$ is determined by the pair of extensions $L$ and $N$. {\em (See \S \ref{sec: u-1 is in u(L)+u(N)} and Theorem \ref{thm: characterization of u-1}.)}
\item[(c)] (1st characterization of $\fu_{-2}(M)$.) Let $\sM^h$ be the second row of \eqref{eq1}, considered as an element of $\Ext^1(\mathbbm{1},\inHom(N,A_1))$. Then $\fu_{-2}(M)$ is the smallest subobject of $\inHom(A_3,A_1)$ such that the pushforward
\[
\sM^h/\fu_{-2}(M) \in \Ext^1(\mathbbm{1}, \inHom(N,A_1)/\fu_{-2}(M) )
\]
of $\sM^h$ along the quotient map $\inHom(N,A_1)\rightarrow \inHom(N,A_1)/\fu_{-2}(M)$ lies in the subgroup $\Ext^1_{\langle L\oplus N\rangle^{\otimes}}(\mathbbm{1}, \inHom(N,A_1)/\fu_{-2}(M) )$ consisting of extensions that belong to the tannakian subcategory $\langle L\oplus N\rangle^{\otimes}$ generated by $L\oplus N$. In particular, $\fu_{-2}(M)$ vanishes if and only if the object $M$ belongs to $\langle L\oplus N\rangle^{\otimes}$. {\em (See Theorem \ref{thm: characterization of u-2} and Corollary \ref{cor: u-2 is zero iff M is generated by L and N}.)}
\item[(d)] (2nd characterization of $\fu_{-2}(M)$, part I - determination of $[\fu(M),\fu(M)]$.) The derived algebra $[\fu(M),\fu(M)]$ of $\fu(M)$ is completely determined by $\fu_{-1}(M)$. More precisely, let $\{ , \}$ be the antisymmetric pairing
\[
V\otimes V \rightarrow \inHom(A_3,A_1)
\]
induced by composition of functions (i.e. given by $\{(f_{12},f_{23}), (g_{12},g_{23})\}=f_{12}g_{23}-g_{12}f_{23}$ after applying a fiber functor $\omega$, where $f_{ij}$ and $g_{ij}$ are linear maps $\omega A_j\rightarrow \omega A_i$). Then $[\fu(M),\fu(M)]=\{\fu_{-1}(M),\fu_{-1}(M)\}$. {\em (See Proposition \ref{prop: [u,u]={u-1,u-1}}.)}
\item[(e)] (Alternative description of $[\fu(M),\fu(M)]$.) The derived algebra $[\fu(M),\fu(M)]$ is the smallest subobject of $\fu_{-2}(M)$ such that the pushforward of the extension \eqref{eq36} along the quotient map $\fu_{-2}(M)\rightarrow \fu_{-2}(M)/[\fu(M),\fu(M)]$ splits. {\em (See Proposition \ref{prop: [u,u]}.)}
\item[(f)] (2nd characterization of $\fu_{-2}(M)$, part II - characterization of $\fu_{-2}(M)/[\fu(M),\fu(M)]$.) There is a canonical isomorphism (independent from any choice of a fiber functor)
\[
\Hom(\frac{\fu_{-2}(M)}{[\fu(M),\fu(M)]}, \inHom(A_3,A_1)) \ \cong \ \Ext^1_{\langle M\rangle^{\otimes}}(A_3,A_1).
\]
In particular, we have $\fu_{-2}(M)=[\fu(M),\fu(M)]$ if and only if there are no nontrivial extensions of $A_3$ by $A_1$ in the subcategory $\langle M\rangle^{\otimes}$. {\em (See Proposition \ref{prop: characterization of u-2/[u,u]} and Corollary \ref{cor: when is [u,u]=u-2}. See also Proposition \ref{prop: map for Ext groups} for the independence of the isomorphism from the choice of a fiber functor.)}
\item[(g)] With the extensions $L$ and $N$ fixed, up to isomorphisms of blended extensions, there exists at most one blended extension $M$ of $N$ by $L$ such that $\fu_{-2}(M)=0$. {\em (See Proposition \ref{prop: at most one M with vanishing u-2})}
\end{itemize}
\end{thm'}

Not every assertion in the theorem requires both hypotheses of semisimplicity of the $A_j$ and \eqref{eq35}. (When the $A_j$ are not semisimple however, one needs to modify the definition of $\fu(M)$; see \S \ref{sec: definition of u}.) See the parts of the paper relevant to each assertion for exactly what assumptions (if any) are needed. In particular, part (b) only needs the hypothesis of semisimplicity of the $A_j$, and part (c) does not need either assumption. Condition \eqref{eq35} is, for instance, automatically satisfied if our category $\bT$ is filtered by weight (as are the categories of motives and mixed Hodge structures) and the $A_j$ are pure objects in an increasing order of weights.\footnote{Importantly, note that unlike the results of \cite{EM2} and \cite{Es23} on motives with maximal unipotent radicals, here one may have
\[\Hom(\inHom(A_3,A_2), \inHom(A_2,A_1))\neq 0.\]}
Note that even if one is only interested in motivic applications, it is preferred to work with conditions like \eqref{eq35} than to assume $\bT$ is filtered by weights, as not every realization category can be equipped with an interesting categorical weight filtration (e.g., there is no such filtration known for the de Rham-Betti realization \cite[\S 7.1.6]{An04}).
\medskip\par 
Part (b) is the most explicit part of Theorem \ref{thm: intro}. Parts (c) and (f), while useful (as illustrated for the latter by the application below), are less explicit. Regarding part (g), one would hope that if $L$ and $N$ are ``blendable" (i.e., if there exists a blended extension of $N$ by $L$) and if $\fu(M)$ is abelian for blended extensions of $N$ by $L$ (a property than depends only on $L$ and $N$, by parts (b) and (d) of the theorem), then there should always exist an $M$ with vanishing $\fu_{-2}(M)$, as shown to be the case by Hardouin \cite{Har05} in the setting of differential operators. In light of Bertrand's work \cite{Ber13} on self-dual blended extensions, one might also expect that for any blendable $L$ and $N$, perhaps under some conditions, there should always exist an $M$ for which $\fu_{-2}(M)=[\fu(M),\fu(M)]$. We have left the investigation of these for the future.

The paper \cite{Es26} is a sequel to this work. Since the original preparation of this manuscript, a generalization of part (b) of Theorem \ref{thm: intro} to filtrations with an arbitrary number of steps and with possibly non-semisimple graded pieces has been given in \cite[Theorem 1.1.1]{Es26} (we caution the reader that in said full generality, the image of the analogue of $(\mathcal{L},\mathcal{N})$ in $\Ext^1(\mathbbm{1}, V/\fu_{-1}(M))$ may not be split; see loc. cit.). Nonetheless, I have kept Theorem \ref{thm: intro}(b) in the present article for the sake of completeness and convenience of the reader, since it is needed for the application discussed below.

\medskip\par 
\noindent {\bf Application to the unipotent part of the Hodge conjecture for 1-motives.} Let $\FF$ be an algebraically closed subfield of $\CC$. The motivic version of the Hodge conjecture for mixed motives predicts that for every (mixed) motive over $\FF$ the motivic Galois group should coincide with the Mumford-Tate group (see \cite[\S 22.3.1]{An04}). Note that for this conjecture to have a precise meaning, a choice of a tannakian category of motives is needed. In \cite{An19} Andr\'{e} has proved the Hodge conjecture for 1-motives, if the motivic Galois group is interpreted in the context of Nori's tannakian category of motives over $\FF$ \cite{HM17}. His proof uses a deformation argument to reduce the problem to the case of semisimple 1-motives. This semisimple case follows by combining the earlier work \cite{An96} of Andr\'{e} in the setting of motives via motivated correspondences with Arapura's result \cite{Ar13} on the equivalence of Andr\'{e}'s category via motivated correspondences and the tannakian subcategory of the category of Nori motives generated by semisimple objects.

Andr\'{e} asks in \cite{An19} whether one can give a proof of the reduction to the semisimple case that does not use deformations. As an application of Theorem \ref{thm: intro} (parts (a), (b), (d), (f)), we propose an approach to such a proof. More precisely, in \S \ref{sec: application to HN for 1-motives}
we prove that the unipotent radicals of the Mumford-Tate and motivic Galois groups of a 1-motive $M$ over $\FF$ coincide, provided that the motivic Galois group is understood in the context of a tannakian category of motives $\mathbf{MM}(\FF)$ over $\FF$ {\it where the Hodge realization map on $\Ext^1_{\langle M\rangle^{\otimes}}(\mathbbm{1}, \QQ(1))$ is injective}. Here, $\langle M\rangle^{\otimes}$ means the tannakian subcategory of $\mathbf{MM}(\FF)$ generated by (the motive of) $M$ (notably, it is not automatic that every extension of $\mathbbm{1}$ by $\QQ(1)$ in $\langle M\rangle^{\otimes}$ comes from a 1-motive, see the next paragraph). Our result for the equality of weight -1 graded parts of the two unipotent radicals is unconditional. The condition of injectivity of the Hodge realization on extension classes of $\mathbbm{1}$ by $\QQ(1)$ in $\langle M\rangle^{\otimes}$ is needed to obtain the equality of weight -2 graded parts. See Theorem \ref{thm: application to Hodge-Nori}.

Deligne \cite[\S 2.4]{De89} conjectures that 1-motives should be closed under extensions in a good category of motives. In particular, every extension of $\mathbbm{1}$ by $\QQ(1)$ in a good category of motives should come from a 1-motive. If this special case of Deligne's conjecture holds true for our category of motives $\mathbf{MM}(\FF)$, then the Hodge realization map is injective on $\Ext^1_{\mathbf{MM}(\FF)}(\mathbbm{1}, \QQ(1))$. Ayoub and Barbieri-Viale have proved in \cite{ABV15} that Deligne's aforementioned conjecture holds in the abelian category of effective Nori motives (see \S 8.11 therein). But unfortunately, one does not know Deligne's conjecture even for the special case of extensions of $\mathbbm{1}$ by $\QQ(1)$ for Nori's tannakian category of (ineffective) motives over $\FF$. 

\medskip\par 
\noindent{\bf A word on conventions and notation.} Our tannakian categories are always assumed to be neutral. Unless there is a chance of misinterpretation, we suppress the category from the notation for Hom and Ext groups. If needed, the category is indicated by a subscript. As already mentioned, $\inHom$ denotes the internal Hom.

The reader can consult \cite[\S 1]{Ber13} or \cite[\S 2]{Es23} for a review of the background material on blended extensions. The collection of isomorphism classes of blended extensions of $N$ by $L$ is denoted by $\Extpan(N,L)$, with the category again dropped from the notation unless there is ambiguity, in which case the intended category will be included as a subscript. With abuse of notation, we use the same symbol for an extension, its middle object, and the corresponding element in the $\Ext^1$ group. Similarly, we use the same notation for a blended extension, it middle object, and the corresponding element in $\Extpan$. If there is a possibility of confusion, we will make the intended interpretation explicit.

Given an extension $\mathcal{E}$ of $\mathbbm{1}$ by $X$ in a tannakian category and a subobject $X'$ of $X$, the notation $\mathcal{E}/X'$ will mean the pushforward of $\mathcal{E}$ along the quotient map $X\rightarrow X/X'$.

Throughout, all actions are designed to be left actions and given an algebraic group $\mathcal{G}$ over a field $\mathbb{K}$, the category of finite-dimensional representations of $\mathcal{G}$ over $\mathbb{K}$ is denoted by $\mathbf{Rep}(\mathcal{G})$. 
\medskip\par 
\noindent{\bf Acknowledgements.} I would like to thank Yves Andr\'{e}, Donu Arapura, Annette Huber, and Kumar Murty for helpful correspondences and conversations. I am also grateful to Cristiana Bertolin, Daniel Bertrand, Charlotte Hardouin, and Peter Jossen as this work was inspired by some of their results. Finally, I would like to express my gratitude to the anonymous referees for their insightful comments, and particularly for suggesting a more natural interpretation of the map of Proposition \ref{prop: map for Ext groups}.

\section{The setup and initial considerations}\label{sec: initial considerations}

\subsection{}\label{sec: data}
From here to the end of \S \ref{sec: characterization of u-2}, we shall fix the following data:
\begin{itemize}
\item[-] A tannakian category $\bT$ over a field $\mathbb{K}$ of characteristic zero. All of what follows takes place in $\bT$.
\item[-] Objects $A_1, A_2$ and $A_3$ of $\bT$. For the time being, there are no conditions on these objects. In particular, they do not have to be semisimple.
\item[-] Two extensions
\[
\begin{tikzcd}[row sep = small]
0 \arrow[r] & A_1 \arrow[r] & L \arrow[r] & A_2 \arrow[r] & 0\\
0 \arrow[r] & A_2 \arrow[r] & N \arrow[r] & A_3 \arrow[r] & 0.
\end{tikzcd}
\]
\item[-] A blended extension $M$ of $N$ by $L$ given by diagram \eqref{eq1}.
\end{itemize}

The notations $\Hom$ and $\Ext$ with no subscripts refer to the Hom and Yoneda Ext groups for $\bT$, and $\inHom$ (also $\inEnd$) refers to internal Homs in $\bT$. Recall that for every $X$ and $Y$ in $\bT$, there is a canonical isomorphism 
\[
\Ext^1(X,Y) \cong \Ext^1(\mathbbm{1}, \inHom(X,Y))
\]
(see for instance, \cite[\S 3.2]{EM1} for an explicit description of this isomorphism). Let
\[
\sL\in \Ext^1(\mathbbm{1},\inHom(A_2,A_1)) \hspace{.3in} (\text{resp.} \ \ \sN\in \Ext^1(\mathbbm{1},\inHom(A_3,A_2)))
\]
be the element corresponding to $L$ (resp. $N$) under the respective canonical isomorphism.

\subsection{}\label{sec: define g(X,Y)}
In this subsection, we recall some well known generalities in tannakian categories. Throughout the paper, unless otherwise indicated by a fiber functor we always mean a fiber functor with values in the category of finite-dimensional vector spaces over $\mathbb{K}$. Let $X$ be an object of $\bT$. Given any fiber functor $\omega$ for $\bT$, we denote the fundamental group of $X$ with respect to $\omega$ by $\mathcal{G}(X,\omega)$; in the standard notation, this is the group scheme $\underline{Aut}^\otimes(\omega |_{\langle X\rangle^\otimes})$ of the tensor automorphisms of the restriction of $\omega$ (and its extensions of scalars) to $\langle X\rangle^\otimes$. The functor $\omega$ gives an equivalence of categories
\begin{equation}\label{eq7}
\langle X\rangle^{\otimes} \rightarrow \mathbf{Rep}(\mathcal{G}(X,\omega))
\end{equation}
(see for instance, \cite{DM82}). Let $\fg(X,\omega)$ be the Lie algebra of $\mathcal{G}(X,\omega)$.

We shall identify $\mathcal{G}(X,\omega)$ as an algebraic subgroup of $\GL(\omega X)$ via the natural embedding
\[
\mathcal{G}(X,\omega) = \underline{Aut}^\otimes(\omega |_{\langle X\rangle^\otimes}) \hookrightarrow \GL(\omega X),
\]
and hence identify the Lie algebra $\fg(X,\omega)$ as a Lie subalgebra of the Lie algebra $\End_\mathbb{K}(\omega X)$ of $\GL(\omega X)$. We identity $\omega\inEnd(M)=\End_\mathbb{K}(\omega X)$, with the action of $\mathcal{G}(X,\omega)$ on $\End_\mathbb{K}(\omega X)$ by conjugation. Via the equivalence of categories \eqref{eq7}, considering the adjoint action of $\mathcal{G}(X,\omega)$ one obtains a canonical Lie subobject 
\[
\fg(X)\subset \inEnd(X)
\]
such that
\[
\omega \fg(X) = \fg(X,\omega)
\]
in $\omega \inEnd(X)=\End_\mathbb{K}(\omega X)$. 

More generally, given any object $Y$ of $\langle X\rangle^{\otimes}$, for any fiber functor $\omega$ set
\[
\mathcal{G}(X,Y,\omega):= \ker\bigm( \mathcal{G}(X,\omega)\twoheadrightarrow \mathcal{G}(Y,\omega)\bigm),
\]
where the surjective arrow is the canonical surjection given by the inclusion of $\langle Y\rangle^{\otimes}$ in $\langle X\rangle^{\otimes}$. Let $\fg(X,Y,\omega)$ be the Lie algebra of $\mathcal{G}(X,Y,\omega)$ (hence identified as a Lie subalgebra of $\End_\mathbb{K}(\omega X)$). Then there exists a canonical Lie subobject
\[
\fg(X,Y)\subset \inEnd(X)
\]
such that for {\it every} fiber functor $\omega$, 
\[
\omega \fg(X,Y) =\fg(X,Y,\omega).
\]
See for example, \S 2.5-2.7 of \cite{EM2} for more details on this, including the independence of $\fg(X,Y)$ (and in particular, $\fg(X)$) from the choice of $\omega$. We have a short exact sequence
\begin{equation}\label{eq22}
\begin{tikzcd}
0 \arrow[r] & \fg(X,Y) \arrow[r] & \fg(X) \arrow[r] & \fg(Y) \arrow[r] & 0  
\end{tikzcd}
\end{equation}
which after applying any fiber functor $\omega$, becomes the exact sequence obtained by applying the Lie algebra functor to the exact sequence
\[
\begin{tikzcd}
1 \arrow[r] & \mathcal{G}(X,Y,\omega) \arrow[r] & \mathcal{G}(X,\omega) \arrow[r] & \mathcal{G}(Y,\omega) \arrow[r] & 1.
\end{tikzcd}
\]

\subsection{}\label{sec: definition of u}
Set
\begin{align*}
\fu(L)&:= \fg(L, A_1\oplus A_2)\\
\fu(N)&:= \fg(N, A_2\oplus A_3)\\
\fu(L\oplus N)&:= \fg(L\oplus N, A_1\oplus A_2\oplus A_3)\\
\fu(M)&:= \fg(M, A_1\oplus A_2\oplus A_3).
\end{align*}

For any choice of fiber functor $\omega$, the image of $\fu(M)$ under $\omega$ is the Lie algebra of the unipotent subgroup
\[
\mathcal{U}(M,\omega):=\mathcal{G}(M, A_1\oplus A_2\oplus A_3, \omega)
\]
of $\mathcal{G}(M,\omega)$. If $A_1, A_2$ and $A_3$ are semisimple, the image of $\fu(M)$ (resp. $\fu(L)$, $\fu(N)$ and $\fu(L\oplus N)$) under $\omega$ is the Lie algebra of the unipotent radical of $\mathcal{G}(M,\omega)$ (resp. $\mathcal{G}(L,\omega)$, $\mathcal{G}(N,\omega)$ and $\mathcal{G}(L\oplus N,\omega)$).

\subsection{}\label{sec: review results for u of extensions}
One has a characterization of $\fu(L)$ as follows. Needless to say, replacing $L$, $A_2$, $A_1$, and $\sL$ respectively by $N$, $A_3$, $A_2$, and $\sN$ we get the statements for $\fu(N)$.

The Lie algebra object $\fu(L)$ is abelian. Indeed, $\fu(L)$ is contained in the abelian Lie subalgebra
\[
\inHom(A_2,A_1)\subset \inEnd(L),
\]
the image of which under every fiber functor $\omega$ is $\Hom_\mathbb{K}(\omega A_2, \omega A_1)$, considered as a subspace of $\End_\mathbb{K}(\omega L)$ via functoriality of $\Hom_\mathbb{K}$. By \cite[Cor. 3.4.1]{EM1} (see also its precursors \cite{Ber01}, \cite{Har06}, \cite{Har11}), we have:

{\it If $A_1$ and $A_2$ are semisimple, then $\fu(L)$ is the smallest subobject of $\inHom(A_2,A_1)$ which satisfies the following two equivalent properties:\\
(1) The extension $\sL$ is the pushforward of an element of $\Ext^1(\mathbbm{1},\fu(L))$ under the inclusion $\fu(L)\subset \inHom(A_2,A_1)$.\\
(2) The extension $\sL/\fu(L)\in \Ext^1(\mathbbm{1}, \inHom(A_2,A_1)/\fu(L))$ ( = pushforward of $\sL$ along the quotient map $\inHom(A_2,A_1)\rightarrow \inHom(A_2,A_1)/\fu(L)$) splits.}

The equivalence of (1) and (2) is seen from the long exact sequence obtained by applying the $\delta$-functor $\Hom(\mathbbm{1}, -)$ to the exact sequence
\[
\begin{tikzcd}
0 \arrow[r] & \fu(L) \arrow[r] & \inHom(A_2,A_1) \arrow[r] & \inHom(A_2,A_1)/\fu(L) \arrow[r] & 0.
\end{tikzcd}
\]

In the more general case where $A_1$ and $A_2$ need not be semisimple one also has a characterization of $\fu(L)$. By \cite[Theorem 3.3.1]{EM1}, we have:

{\it The object $\fu(L)$ is the smallest subobject of $\inHom(A_2,A_1)$ which satisfies the following property: the pushforward $\sL/\fu(L)\in \Ext^1(\mathbbm{1}, \inHom(A_2,A_1)/\fu(L))$ belongs to the subgroup 
\[
\Ext^1_{\langle A_1,A_2\rangle^{\otimes}}(\mathbbm{1}, \inHom(A_2,A_1)/\fu(L)) \subset \Ext^1(\mathbbm{1}, \inHom(A_2,A_1)/\fu(L))
\]
(where $\langle A_1,A_2\rangle^{\otimes}$ is the tannakian subcategory of $\bT$ generated by $A_1$ and $A_2$).}

\subsection{}
We shall consider the injective arrows in \eqref{eq1} as inclusion maps. Our blended extension gives a 3-step filtration
\[
W_{-3}M:= 0 \subset W_{-2}M:= A_1 \subset W_{-1}M:= L \subset W_{0}M:=M
\]
of $M$ and a graded isomorphism
\begin{equation}\label{eq3}
\Gr(M) \cong A_1\oplus A_2\oplus A_3,
\end{equation}
where $\Gr(M)$ is the associated graded of $M$ with respect to $W_\bullet$. 

Let
\[
W_{-1}\inEnd(M) \subset \inEnd(M)
\]
be the subobject whose image under any fiber functor $\omega$ consists of the maps $f\in \End_\mathbb{K}(\omega M)$ such that $f(\omega W_nM) \subset \omega W_{n-1}M$ for each $n$, i.e.
\[
f(\omega M) \subset \omega L, \ \ \ f(\omega L)\subset \omega A_1, \ \ \ f(\omega A_1) = 0.
\]
Then
\[
\fu(M) \subset W_{-1}\inEnd(M).
\]
Indeed, for every fiber functor $\omega$, every element of the fundamental group $\mathcal{G}(M,\omega)$ preserves the filtration $\omega W_\bullet M$ of $\omega M$, so that $\mathcal{G}(M,\omega)$ is contained in the parabolic subgroup of $\GL(\omega M)$ corresponding to this filtration. The subgroup $\mathcal{U}(M,\omega)$ of $\mathcal{G}(M,\omega)$ (i.e., the kernel of the natural surjection of $\mathcal{G}(M,\omega)$ onto $\mathcal{G}(A_1\oplus A_2\oplus A_3,\omega)$) is contained in the unipotent radical of this parabolic subgroup, consisting of those elements of the parabolic subgroup that induce identity on $\Gr(\omega M)$. It follows that $\omega \fu(M)$ is contained in the Lie algebra of this unipotent radical, which is exactly $\omega W_{-1}\inEnd(M)$. 

\subsection{}\label{sec: construction of pi}
The bifunctoriality of $\inHom$ gives a canonical embedding
\[
\inHom(A_3,A_1) \subset W_{-1}\inEnd(M).
\]
Applying any fiber functor $\omega$, this inclusion sends a linear map $\omega A_3\rightarrow \omega A_1$ to the composition
\[
\omega M\twoheadrightarrow \omega A_3\rightarrow \omega A_1 \hookrightarrow \omega M.
\]
This embedding identifies $\omega \inHom(A_3,A_1)$ as the subspace of $\omega W_{-1}\inEnd(M)$ consisting of linear maps $\omega M\rightarrow \omega M$ which vanish on $\omega L$ and whose image is in $\omega A_1$. It fits in a canonical short exact sequence
\begin{equation}\label{eq37}
\begin{tikzcd}
0 \arrow[r] & \inHom(A_3,A_1) \arrow[r] & W_{-1}\inEnd(M) \arrow[r, "\pi"] & \inHom(A_2,A_1)\oplus \inHom(A_3,A_2) \arrow[r] & 0,
\end{tikzcd}
\end{equation}
where the map 
\[
\pi : W_{-1}\inEnd(M) \twoheadrightarrow \inHom(A_2,A_1)\oplus \inHom(A_3,A_2)
\]
is obtained as follows: Bifunctoriality of $\inHom$ gives rise to two exact sequences
\[
\begin{tikzcd}
0 \arrow[r] & \inHom(A_3,A_1) \arrow[r] & \inHom(N,A_1) \arrow[r, "\pi_{1}"] & \inHom(A_2,A_1) \arrow[r] & 0
\end{tikzcd}
\]
and
\[
\begin{tikzcd}
0 \arrow[r] & \inHom(A_3,A_1) \arrow[r] & \inHom(A_3,L) \arrow[r, "\pi_2"] & \inHom(A_3,A_2) \arrow[r] & 0.
\end{tikzcd}
\]
We also have embeddings
\begin{align*}
\inHom(N,A_1) &\subset W_{-1}\inEnd(M)\\
\inHom(A_3,L) &\subset W_{-1}\inEnd(M)
\end{align*}
the sum of which gives rise to a canonical isomorphism
\[
\frac{\inHom(N,A_1)\oplus \inHom(A_3,L)}{\text{anti-diagonal copy of} \ \inHom(A_3,A_1)} \ \cong \ W_{-1}\inEnd(M).
\]
Thus $W_{-1}\inEnd(M)$ is the fibered coproduct of $\inHom(N,A_1)$ and $\inHom(A_3,L)$ over $\inHom(A_3,A_1)$. The map $\pi$ is now obtained from the compositions
\[
\inHom(N,A_1) \stackrel{\pi_1}{\twoheadrightarrow} \inHom(A_2,A_1) \hookrightarrow \inHom(A_2,A_1)\oplus \inHom(A_3,A_2)
\]
and
\[
\inHom(A_3,L) \stackrel{\pi_2}{\twoheadrightarrow} \inHom(A_3,A_2) \hookrightarrow \inHom(A_2,A_1)\oplus \inHom(A_3,A_2)
\]
(the kernel of each of which is $\inHom(A_3,A_1)$) via the universal property of a fibered coproduct.

After applying a fiber functor $\omega$, the map $\pi$ has the following simple description: given any linear map $f\in W_{-1}\End_\mathbb{K}(\omega M)$, write $f$ as $f_1+f_2$ for some 
\[
f_1\in \Hom_\mathbb{K}(\omega N, \omega A_1)\subset W_{-1}\End_\mathbb{K}(\omega M)
\] 
and
\[
f_2\in \Hom_\mathbb{K}(\omega A_3, \omega L)\subset W_{-1}\End_\mathbb{K}(\omega M).
\]
The elements $f_1$ and $f_2$ are not unique but their images under $\pi_1$ and $\pi_2$ are unique. We have $\pi(f)=(\pi_1(f_1),\pi_2(f_2))$. 

Let 
\[
\pi_{12}: W_{-1}\inEnd(M) \rightarrow \inHom(A_2,A_1) \ \ \ \text{and} \ \ \ \pi_{23}: W_{-1}\inEnd(M) \rightarrow \inHom(A_3,A_2)
\]
be the compositions of $\pi$ with the projections maps to the two factors of $\inHom(A_2,A_1)\oplus \inHom(A_3,A_2)$.

\subsection{}\label{sec: def of u-1 and u-2}
Set 
\begin{align*}
\fu_{-2}(M) &:= \fu(M)\cap \inHom(A_3,A_1)\\
\fu_{-1}(M) &:= \pi(\fu(M)) \subset \inHom(A_2,A_1)\oplus \inHom(A_3,A_2).
\end{align*}
We have then a commutative diagram
\begin{equation}\label{eq10}
\begin{tikzcd}
0 \arrow[r] & \fu_{-2}(M) \arrow[d, hookrightarrow] \arrow[r] & \fu(M) \arrow[d, hookrightarrow] \arrow[r, "\pi |_{\fu(M)}"] & \fu_{-1}(M) \arrow[d, hookrightarrow] \arrow[r] & 0\\
0 \arrow[r] & \inHom(A_3,A_1) \arrow[r] & W_{-1}\inEnd(M) \arrow[r, "\pi"] & \inHom(A_2,A_1)\oplus \inHom(A_3,A_2) \arrow[r] & 0
\end{tikzcd}
\end{equation}
with exact rows and inclusion vertical arrows.

\subsection{}
In general, the two subobjects $\fu_{-2}(M)$ and $\fu_{-1}(M)$ of $\inHom(A_3,A_1)$ and $\inHom(A_2,A_1)\oplus \inHom(A_3,A_2)$ may not determine the subobject $\fu(M)$ of $W_{-1}\inEnd(M)$. They will however, under a mild condition that often holds in practical situations of interest.

\begin{prop}\label{prop: when u is determined by its two subquotients}
Let $\omega$ be a fiber functor. Set $\mathcal{G}:=\mathcal{G}(M,\omega)$. Suppose that there exists a reductive subgroup $\mathcal{R}$ of $\mathcal{G}$ such that 
\begin{equation}\label{eq2}
\Hom_{\mathcal{R}}(\omega \fu_{-1}(M), \omega \inHom(A_3,A_1)) \ = \ 0
\end{equation}
(equivalently, such that $\omega \fu_{-1}(M)$ and $\omega \inHom(A_3,A_1)$, when considered as representations of $\mathcal{R}$, do not have any isomorphic subquotients). Then $\fu(M)$ is uniquely determined inside $W_{-1}\inEnd(M)$ by the two subobjects $\fu_{-2}(M)$ and $\fu_{-1}(M)$ of $\inHom(A_3,A_1)$ and $\inHom(A_2,A_1)\oplus \inHom(A_3,A_2)$.
\end{prop}

\begin{proof}
We will argue that $\omega \fu(M)$ is determined by $\omega \fu_{-2}(M)$ and $\omega \fu_{-1}(M)$. We have a diagram
\[
\begin{tikzcd}
0 \arrow[r] & \fu_{-2}(M) \arrow[d, hookrightarrow] \arrow[r] & \fu(M) \arrow[d, hookrightarrow] \arrow[r, "\pi |_{\fu(M)}"] & \fu_{-1}(M) \ar[equal]{d} \arrow[r] & 0\\
0 \arrow[r] & \inHom(A_3,A_1) \ar[equal]{d} \arrow[r] & \pi^{-1}\fu_{-1}(M) \arrow[d, hookrightarrow] \arrow[r, "\pi |_{\pi^{-1}\fu_{-1}(M)}"] & \fu_{-1}(M) \arrow[d, hookrightarrow] \arrow[r] & 0\\
0 \arrow[r] & \inHom(A_3,A_1) \arrow[r] & W_{-1}\inEnd(M) \arrow[r, "\pi"] & \inHom(A_2,A_1)\oplus \inHom(A_3,A_2) \arrow[r] & 0,
\end{tikzcd}
\]
where the second row is the pullback of the third row along the inclusion of $\fu_{-1}M$. Given any section $s$ of $\omega \pi |_{\pi^{-1}\fu_{-1}(M)}$ in the category of vector spaces that maps $\omega \fu_{-1}(M)$ into $\omega \fu(M)$, we have $\omega \fu(M)=\omega\fu_{-2}(M)+s(\omega \fu_{-1}(M))$ in $\omega W_{-1}\inEnd(M)$. Since $\mathcal{R}$ is reductive, the map $\omega \pi |_{\pi^{-1}\fu_{-1}(M)}$ admits an $\mathcal{R}$-equivariant section $s_0$, which is unique thanks to \eqref{eq2}. We will argue that $s_0$ maps $\omega \fu_{-1}(M)$ into $\omega \fu(M)$, hence establishing the result. Indeed, since $\mathcal{R}$ is reductive, the map $\omega \pi |_{\fu(M)}$ admits an $\mathcal{R}$-equivariant section, which after compositing with the inclusion $\omega \fu(M)\subset \omega \pi^{-1}\fu_{-1}(M)$ gives an $\mathcal{R}$-equivariant section of $\omega \pi |_{\pi^{-1}\fu_{-1}(M)}$. This section must coincide with $s_0$ by the uniqueness property. 
\end{proof}

\begin{rem}
The condition of Proposition \ref{prop: when u is determined by its two subquotients} for $\fu(M)$ to be determined by $\fu_{-2}(M)$ and $\fu_{-1}(M)$ holds for instance, if $\mathbf{T}$ has a functorial weight filtration (i.e., a filtration similar to the weight filtration in the category of mixed Hodge structures: indexed by $\ZZ$, finite increasing on every object, functorial, exact, and compatible with tensor products), and $A_1,A_2,A_3$ are pure and of increasing order of weights. The condition then holds for $\mathcal{R}$ taken to be the image of the weight cocharacter $\mathbb{G}_m\rightarrow \mathcal{G}(M,\omega)$.
\end{rem}

\subsection{}\label{sec: splitting} 
Let $\omega$ be a fiber functor. Consider the filtration
\[
0\subset \omega A_1\subset \omega A_1\oplus \omega A_2 \subset \omega A_1\oplus \omega A_2\oplus \omega A_3
\]
of $\omega A_1\oplus \omega A_2\oplus \omega A_3$. Applying $\omega$ to our blended extension \eqref{eq1}, we obtain a filtration
\[
0\subset \omega A_1\subset \omega L\subset \omega M
\] 
of $\omega M$ and a graded isomorphism
\[
\Gr(\omega M) \cong \omega A_1\oplus \omega A_2\oplus \omega A_3.\]
By a {\it splitting} of $\omega M$ we mean an isomorphism of filtered vector spaces
\[
\varphi: \omega M \rightarrow \omega A_1\oplus \omega A_2\oplus \omega A_3,
\]
such that the induced isomorphism 
\[\Gr(\varphi):\Gr(\omega M)\rightarrow \Gr(\omega A_1\oplus \omega A_2\oplus \omega A_3) \cong \omega A_1\oplus \omega A_2\oplus \omega A_3\]
coincides with the isomorphism coming from the blended extension \eqref{eq1}. The data of a splitting of $\omega M$ is equivalent to the data of a section of each of the structure maps $\omega L\twoheadrightarrow \omega A_2$ and $\omega M\twoheadrightarrow \omega A_3$. In particular, a splitting of $\omega M$ induces isomorphisms of vector spaces $\omega L\cong \omega A_1\oplus \omega A_2$ and $\omega N\cong \omega A_2\oplus \omega A_3$.

\subsection{}\label{sec: description of maps wrt a splitting} Let $\omega$ be a fiber functor. Choose a splitting $\varphi: \omega M\xrightarrow{\ \simeq \ } \omega A_1\oplus \omega A_2\oplus \omega A_3$. Using $\varphi$ to identity $\omega M$ with $\omega A_1\oplus \omega A_2\oplus \omega A_3$, we can write elements of 
\[\End_\mathbb{K}(\omega M) = \End_\mathbb{K}(\omega A_1\oplus \omega A_2\oplus \omega A_3) \cong \bigoplus_{1\leq i,j\leq 3} \Hom_\mathbb{K}(\omega A_j, \omega A_i)\]
as $3\times 3$ matrices $f=(f_{ij})$, where $f_{ij}$ is the component of $f$ in $\Hom_\mathbb{K}(\omega A_j, \omega A_i)$. The Lie subalgebra $\omega W_{-1}\inEnd(M)$ is then the space of all strictly upper triangular elements of $\End_\mathbb{K}(\omega M)$. The Lie bracket of $\omega W_{-1}\inEnd(M)$ is simply the usual Lie bracket on matrices. 

The canonical embedding
\[
\inHom(A_3, A_1) \hookrightarrow W_{-1}\inEnd(M)
\]
after applying $\omega$ simply places $f_{13}: \omega A_3\rightarrow \omega A_1$ as the (13)-entry of a matrix. The canonical surjection 
\[
\pi=(\pi_{12}, \pi_{23}): W_{-1}\inEnd(M) \twoheadrightarrow \inHom(A_2,A_1)\oplus \inHom(A_3,A_2)
\]
after applying $\omega$ is then simply the map
\[
(f_{ij}) \mapsto (f_{12},f_{23}).
\]
The two maps $\omega \pi_{12}$ and $\omega \pi_{23}$ simply send $(f_{ij})$ to $f_{12}$ and $f_{23}$, respectively. Note that a different choice of splitting would result to conjugation of $(f_{ij})$ by a unipotent upper triangular matrix in $\GL(\omega A_1\oplus \omega A_2\oplus \omega A_3)$, and hence indeed the (12) and (23) entries would not change.

Use $\varphi$ to also identity $\omega L$ and $\omega N$ respectively with $\omega A_1\oplus \omega A_2$ and $\omega A_2\oplus \omega A_3$. Identifying $\mathcal{G}(L,\omega)$, $\mathcal{G}(M,\omega)$ and $\mathcal{G}(N,\omega)$ as subgroups of 
\[
\GL(\omega L)\cong \GL(\omega A_1\oplus \omega A_2),  \ \GL(\omega M)\cong \GL(\omega A_1\oplus \omega A_2\oplus \omega A_3) \ \ \text{and} \ \GL(\omega N)\cong \GL(\omega A_2\oplus \omega A_3),
\]
the canonical surjections
\[
\mathcal{G}(M,\omega)\twoheadrightarrow \mathcal{G}(L,\omega) \ \ \ \text{and} \ \ \ \mathcal{G}(M,\omega)\twoheadrightarrow \mathcal{G}(N,\omega)
\]
arising from the inclusions of $\langle L \rangle^{\otimes}$ and $\langle N \rangle^{\otimes}$ in $\langle M \rangle^{\otimes}$ are respectively simply given by
\[
\begin{pmatrix}
\sigma_{11}& \sigma_{12}&\sigma_{13}\\
0 & \sigma_{22}&\sigma_{23}\\
0&0&\sigma_{33}\end{pmatrix}\mapsto \begin{pmatrix}
\sigma_{11}& \sigma_{12}\\
0 & \sigma_{22}\end{pmatrix} \ \ \ \text{and} \ \ \ \begin{pmatrix}
\sigma_{11}& \sigma_{12}&\sigma_{13}\\
0 & \sigma_{22}&\sigma_{23}\\
0&0&\sigma_{33}\end{pmatrix}\mapsto \begin{pmatrix}
\sigma_{22}& \sigma_{23}\\
0 & \sigma_{33}\end{pmatrix}.
\]

\subsection{}\label{sec: u-1 is in u(L)+u(N)}
It is not difficult to see that the subobject $\fu_{-1}(M)$ of $\inHom(A_2,A_1)\oplus \inHom(A_3,A_2)$ only depends on $L$ and $N$, and not on the choice of the blended extension $M$. Indeed, fix a fiber functor $\omega$. Dropping $\omega$ from the notation for tannakian groups for simplicity, we have a diagram
\[
\begin{tikzcd}
\mathcal{G}(M) \arrow[r, twoheadrightarrow] \arrow[d, twoheadrightarrow] & \mathcal{G}(L\oplus N) \arrow[r, hookrightarrow] \arrow[d, twoheadrightarrow] & \mathcal{G}(L)\times \mathcal{G}(N) \arrow[d, twoheadrightarrow]\\
\mathcal{G}(A_1\oplus A_2\oplus A_3) \ar[r, equal] & \mathcal{G}(A_1\oplus A_2\oplus A_3) \arrow[r, hookrightarrow] & \mathcal{G}(A_1\oplus A_2)\times \mathcal{G}(A_2\oplus A_3),
\end{tikzcd}
\]
where all the arrows are induced by restrictions to the respective subcategories. Focusing on the induced maps between the kernels of the vertical maps and passing to the Lie algebras, we obtain $\mathcal{G}(M)$-equivariant maps
\[
\omega \fu(M) \twoheadrightarrow \omega \fu(L\oplus N) \hookrightarrow \omega \fu(L)\oplus \omega \fu(N)
\]
and hence morphisms
\begin{equation}\label{eq24}
\fu(M) \twoheadrightarrow \fu(L\oplus N) \hookrightarrow \fu(L)\oplus \fu(N),
\end{equation}
which are in fact, independent of the choice of $\omega$. There is a commutative diagram
\begin{equation}\label{eq23}
\begin{tikzcd}
\fu(M) \arrow[r, twoheadrightarrow] \arrow[d, hookrightarrow] & \fu(L\oplus N) \arrow[r, hookrightarrow] & \fu(L)\oplus \fu(N) \arrow[d, hookrightarrow] \\
W_{-1}\inEnd(M) \arrow[rr, "\text{$\pi=(\pi_{12},\pi_{23})$}"] && \inHom(A_2,A_1)\oplus \inHom(A_3,A_2).
\end{tikzcd}
\end{equation}
That this diagram commutes is seen immediately upon applying $\omega$ and taking a splitting of $\omega M$. We thus have 
\begin{equation}\label{eq25}
\fu_{-1}(M) = \text{Im} \bigm(\fu(L\oplus N) \hookrightarrow \fu(L)\oplus \fu(N) \bigm).
\end{equation}
In particular, we obtain:
\begin{prop}
The subobject $\fu_{-1}(M)$ of $\inHom(A_2,A_1)\oplus \inHom(A_3,A_2)$ only depends on the top horizontal and right vertical extensions $L$ and $N$ of \eqref{eq1}. (That is, it does not depend on the choice of $M$ in $\Extpan(N,L)$).
\end{prop}

\begin{notation}\label{notation piL and piN}
Denote the compositions of \eqref{eq24} with the projections onto $\fu(L)$ and $\fu(N)$ respectively by $\pi_L$ and $\pi_N$. The maps $\pi_L$ and $\pi_N$ are respectively induced by the inclusions of $\langle L\rangle^{\otimes}$ and $\langle N\rangle^{\otimes}$ in $\langle M\rangle^{\otimes}$.
\end{notation}

\subsection{}\label{sec: R-equivariant splitting}
Let $\omega$ be a fiber functor for $\bT$ and $\mathcal{R}$ a subgroup of $\mathcal{G}(M,\omega)$. A splitting $\omega M\rightarrow \omega A_1\oplus\omega A_2\oplus \omega A_3$ is said to be $\mathcal{R}$-equivariant if it is $\mathcal{R}$-equivariant with respect to the natural action of $\mathcal{R}$ on $\omega M$ (restricted from the action of $\mathcal{G}(M,\omega)$) and the direct sum action of $\mathcal{R}$ on $\omega A_1\oplus\omega A_2\oplus \omega A_3$. If a splitting is $\mathcal{R}$-equivariant, then the induced isomorphisms $\omega L\cong \omega A_1\oplus\omega A_2$ and $\omega N\cong \omega A_2\oplus\omega A_3$ are also $\mathcal{R}$-equivariant with respect to the obvious actions. An $\mathcal{R}$-equivariant splitting of $\omega M$ exists if and only if the surjections $\omega L\twoheadrightarrow \omega A_2$ and $\omega M\twoheadrightarrow \omega A_3$ admit $\mathcal{R}$-equivariant sections. In particular, if $\mathcal{R}$ is reductive, then there always exists an $\mathcal{R}$-equivariant splitting of $\omega M$.

\begin{lemma}\label{lem: the obvious R-equiv section of pi}
Let $\mathcal{R}$ be a subgroup of $\mathcal{G}(M,\omega)$ such that an $\mathcal{R}$-equivariant splitting of $\omega M$ exists. Identifying $\omega M$ with $\omega A_1\oplus \omega A_2\oplus \omega A_3$ via an $\mathcal{R}$-equivariant splitting $\varphi$, the map
\[
\omega \inHom(A_2,A_1) \oplus \omega \inHom(A_3,A_2) \rightarrow \omega W_{-1}\inEnd(M) \hspace{.3in} (f_{12},f_{23})\mapsto \begin{pmatrix}0 & f_{12} &0 \\ 0& 0& f_{23} \\
0 & 0 &0\end{pmatrix}
\]
is $\mathcal{R}$-equivariant.
\end{lemma}
\begin{proof}
Let $\sigma\in \mathcal{G}(M,\omega)$. We write $\sigma$ as a $3\times 3$ matrix via the identifications $\mathcal{G}(M,\omega)\subset \GL(\omega M)$ and $\omega M\cong \omega A_1\oplus \omega A_2\oplus \omega A_3$ given by the splitting $\varphi$. Then $\sigma=(\sigma_{ij})$ is upper triangular, and its diagonal entry $\sigma_{jj}$ is the natural action of $\sigma$, as an automorphism of $\omega$, on $\omega A_j$. The actions of $\mathcal{G}(M,\omega)$ corresponding to $\inHom(A_2,A_1)$, $\inHom(A_3,A_2)$ and $W_{-1}\inEnd(M)$ are given as follows: for any linear maps $f_{12}:\omega A_2\rightarrow \omega A_1$ and $f_{23}:\omega A_3\rightarrow \omega A_2$, we have $\sigma\cdot f_{12}=\sigma_{11}f_{12}\sigma_{22}^{-1}$ and $\sigma\cdot f_{23}=\sigma_{22}f_{23}\sigma_{33}^{-1}$. For any $f\in \omega W_{-1}\inEnd(M)$, we have $\sigma\cdot f=(\sigma_{ij})f(\sigma_{ij})^{-1}$. The $\mathcal{R}$-equivariance of the map in the statement of the lemma is checked by a direct computation on noting that $(\sigma_{ij})$ is diagonal if $\sigma\in\mathcal{R}$.
\end{proof}

\subsection{}
Combining Propositions \ref{prop: when u is determined by its two subquotients} (and its proof) with Lemma \ref{lem: the obvious R-equiv section of pi} we obtain the following:

\begin{prop}\label{prop: decomposition of u as superdiagonals plus top right entries}
Let $\omega$ be a fiber functor for $\bT$. Set $\mathcal{G}:=\mathcal{G}(M,\omega)$. Suppose that there exists a reductive subgroup $\mathcal{R}$ of $\mathcal{G}$ such that condition \eqref{eq2} of Proposition \ref{prop: when u is determined by its two subquotients} holds. Identifying $\omega M$ with $\omega A_1\oplus \omega A_2\oplus \omega A_3$ via an $\mathcal{R}$-equivariant splitting $\omega M\rightarrow \omega A_1\oplus \omega A_2\oplus \omega A_3$, for every (strictly upper triangular) $f=(f_{ij})\in \omega W_{-1}\inEnd(M)$ we have
\[
f\in \omega\fu(M) \ \ \ \ \Longleftrightarrow \ \ \ \begin{pmatrix}
0 & f_{12} & 0\\
0 & 0 & f_{23}\\
0&0&0\end{pmatrix} \ \text{and} \ \begin{pmatrix}
0 & 0 & f_{13}\\
0 & 0 & 0\\
0&0&0\end{pmatrix}  \in  \omega\fu(M).
\]
\end{prop}

\begin{proof}
Recall from the proof of Proposition \ref{prop: when u is determined by its two subquotients} that we have an internal direct sum decomposition of vector spaces
\[
\omega \fu(M) = \omega\fu_{-2}(M)+s(\omega\fu_{-1}(M)) 
\]
in $\omega W_{-1}\inEnd(M)$, where $s$ is the unique $\mathcal{R}$-equivariant section of $\omega \pi|_{\pi^{-1}\fu_{-1}(M)}$ (notation as in the proof of Proposition \ref{prop: when u is determined by its two subquotients}). By the uniqueness of this section, $s$ must be the restriction to $\omega \fu_{-1}(M)$ of the section of $\pi$ constructed in Lemma \ref{lem: the obvious R-equiv section of pi}. The result follows.
\end{proof}

\begin{rem}
The reader should be alert that neither the result above nor Proposition \ref{prop: when u is determined by its two subquotients} asserts that (even under the hypotheses of the results) $\fu(M)$ decomposes as a direct sum of $\fu_{-2}(M)$ and $\fu_{-1}(M)$ in $\bT$.
\end{rem}

\subsection{}
For every object $X$ of $\bT$ the object $\fg(X)$ is a Lie subobject of $\inEnd(X)$. That is, one has a diagram
\[
\begin{tikzcd}
\fg(X)\otimes \fg(X) \arrow[d, hookrightarrow] \arrow[r, "\text{$[ \ , \ ]$}" ] & \fg(X) \arrow[d, hookrightarrow]\\
\inEnd(X) \otimes \inEnd(X) \arrow[r, "\text{$[ \ , \ ]$}"] & \inEnd(X)
\end{tikzcd}
\]
in $\bT$ whose image under every $\omega$ is a similar diagram for the Lie brackets in the classical sense of the Lie algebras $\fg(X,\omega)$ and $\End_\mathbb{K}(\omega X)$, making the former a Lie subalgebra of the latter. The Lie bracket
\[
\End_\mathbb{K}(\omega X)\otimes \End_\mathbb{K}(\omega X) \rightarrow \End_\mathbb{K}(\omega X)
\]
is simply the usual Lie bracket of endomorphisms, given by 
\begin{equation}\label{eq9}
[f,f'] = f\circ f'-f'\circ f.
\end{equation}

In the case of our blended extension $M$, the above diagram of Lie brackets restricts to a diagram
\[
\begin{tikzcd}
\fu(M)\otimes \fu(M) \arrow[d, hookrightarrow] \arrow[r, "\text{$[ \ , \ ]$}" ] & \fu(M) \arrow[d, hookrightarrow]\\
W_{-1}\inEnd(M) \otimes W_{-1}\inEnd(M) \arrow[r, "\text{$[ \ , \ ]$}"] & W_{-1}\inEnd(M).
\end{tikzcd}
\]
By \eqref{eq9}, the derived Lie algebra $[W_{-1}\inEnd(M), W_{-1}\inEnd(M)]$ is contained in the subobject $\inHom(A_3,A_1)$ of $W_{-1}\inEnd(M)$ (and they are equal if $A_2\neq 0$). Thus
\[
[\fu(M),\fu(M)] \subset \inHom(A_3,A_1) \cap \fu(M) = \fu_{-2}(M).
\]

\subsection{}\label{sec:summary of setup}
Let us summarize our picture so far. Given our fixed data of \S \ref{sec: data}, we have
\[
0\subset [\fu(M),\fu(M)] \subset \fu_{-2}(M) = \fu(M)\cap \inHom(A_3,A_1)\subset \fu(M) \subset W_{-1}\inEnd(M).
\]
We also have
\[
\frac{\fu(M)}{\fu_{-2}(M)}\stackrel{\pi}{\cong} \fu_{-1}(M):=\pi(\fu(M))\subset \inHom(A_2,A_1)\oplus \inHom(A_3,A_2),
\]
where
\[
\pi: W_{-1}\inEnd(M) \twoheadrightarrow \inHom(A_2,A_1)\oplus \inHom(A_3,A_2)
\]
is the canonical surjection constructed in \S \ref{sec: construction of pi} (with the simple description $(f_{ij})\mapsto (f_{12},f_{23})$ after applying a fiber functor and taking a splitting). Moreover, we have seen that $\fu_{-1}(M)$ only depends of $L$ and $N$, and is contained in the subobject $\fu(L)\oplus \fu(N)$ of $\inHom(A_2,A_1)\oplus \inHom(A_3,A_2)$.

When the mild condition given in Proposition \ref{prop: when u is determined by its two subquotients} holds for a fiber functor $\omega$, then the subobjects $\fu_{-2}(M)$ and $\fu_{-1}(M)$ of $\inHom(A_3,A_2)$ and $\inHom(A_2,A_1)\oplus \inHom(A_3,A_2)$ uniquely determine the subobject $\fu(M)$ of $W_{-1}\inEnd(M)$. In this situation, for a suitable choice of splitting $\omega M\rightarrow \omega A_1\oplus \omega A_2\oplus \omega A_3$, we also had a concrete decomposition of $\omega \fu(M)$ as a vector space in Proposition \ref{prop: decomposition of u as superdiagonals plus top right entries}.

\subsection{}\label{sec: conditions C1 and C2} 
Our objective in this paper is to study 
\begin{center}$\fu_{-1}(M)$, $\fu_{-2}(M)$, $[\fu(M),\fu(M)]$ and $\fu_{-2}(M)/[\fu(M),\fu(M)]$.\end{center}
We shall mostly (to be made clear as we go through the paper) work under the following assumption:
\begin{itemize}
\item[(C1)] The objects $A_1, A_2$ and $A_3$ are semisimple.
\end{itemize}
This condition is equivalent to $\mathcal{G}(A_1\oplus A_2\oplus A_3,\omega)$ being reductive for any fiber functor $\omega$, and also is equivalent to $\mathcal{U}(M,\omega)$ being the unipotent radical of $\mathcal{G}(M,\omega)$ for any such $\omega$.
\medskip\par 
In addition, in our study of $\fu_{-2}(M)/[\fu(M),\fu(M)]$ (\S \ref{sec: char of u-2/[u,u]} and \S \ref{sec: at most one M with vanishing u-2}) we shall also assume the following condition:
\begin{itemize}
\item[(C2)] We have 
\[
\Hom(A_3\otimes A_2, A_3\otimes A_1) \cong \Hom(A_1\otimes A_3, A_1\otimes A_2) \cong 0.
\]
\end{itemize}
In view of the canonical isomorphisms $\Hom(X,Y)\cong \Hom(\mathbbm{1},\inHom(X,Y))$ and $\inHom(X,Y)\cong X^\vee\otimes Y$ in a tannakian category, (C2) is equivalent to
\begin{equation}\label{eq16}
\Hom(\inHom(A_2,A_1), \inHom(A_3,A_1)) \cong \Hom(\inHom(A_3,A_2), \inHom(A_3,A_1)) \cong 0.
\end{equation}

Both (C1) and (C2) are satisfied in a significant class of interesting practical situations, e.g., if $\bT$ has a functorial weight filtration and $A_1, A_2, A_3$ are semisimple and pure of increasing order of weights. We highlight that under (C1) and (C2) we may still have
\[
\Hom(\inHom(A_3,A_2), \inHom(A_2,A_1)) \neq 0.
\]
(Compare with \eqref{eq16}.) If we were to assume that the last Hom group is zero, then the study of $\fu_{-1}(M)$ would considerably simplify, as it would bring us to a situation similar to the ``graded-independent" case in \cite{EM2} and \cite{Es23}. In this simplified situation, we will always have $\fu_{-1}(M)=\fu(L)\oplus \fu(N)$, as can be seen from Theorem \ref{thm: characterization of u-1}. One of the main goals of this paper is to go beyond this simplified and limiting case.

Conditions (C1) and (C2) together guarantee that the condition of Proposition \ref{prop: when u is determined by its two subquotients} holds and hence $\fu(M)$ is completely determined by $\fu_{-1}(M)$ and $\fu_{-2}(M)$. Indeed, referring to the notation of Proposition \ref{prop: when u is determined by its two subquotients}, take $\mathcal{R}$ to be a Levi factor of $\mathcal{G}(M,\omega)$. Then we have
\[
\Hom_\mathcal{R}(\omega \inHom(A_j,A_{j-1}), \omega\inHom(A_3,A_1))\cong \Hom(\inHom(A_j,A_{j-1}), \inHom(A_3,A_1)) \cong 0
\]
for $j=2,3$.

\subsection{} We end this section with a simple observation. Suppose $A_2\neq 0$. Let $\omega$ be a fiber functor. Then the only Lie subalgebra of $W_{-1}\End_\mathbb{K}(\omega M)$ that surjects onto $\Hom_\mathbb{K}(\omega A_2,\omega A_1)\oplus \Hom_\mathbb{K}(\omega A_3,\omega A_2)$ by $\omega \pi$ is $W_{-1}\End_\mathbb{K}(\omega M)$. Indeed, let $\fv$ be a Lie subalgebra of $W_{-1}\End_\mathbb{K}(\omega M)$ whose image under $\omega\pi$ contains $(f_{12},g_{23})$ for every $f_{12}:\omega A_2\rightarrow \omega A_1$ and $g_{23}:\omega A_3\rightarrow \omega A_2$. Choose a splitting of $\omega M$ to write elements of $W_{-1}\End_\mathbb{K}(\omega M)$ as strictly upper triangular matrices. By the description of $\omega \pi$ given in \S \ref{sec: description of maps wrt a splitting} and surjectivity of $\omega\pi$ on $\fv$, for every $f_{12}:\omega A_2\rightarrow \omega A_1$ and $g_{23}:\omega A_3\rightarrow \omega A_2$ there exist $f_{13}, g_{13}\in \Hom_\mathbb{K}(\omega A_3, \omega A_1)$ such that $\fv$ contains
\[
f:= \begin{pmatrix}
0 & f_{12} &  f_{13}\\
0 & 0 &0\\
0&0&0\end{pmatrix} \ \ \text{and} \ \ g:= \begin{pmatrix}
0 & 0 &  g_{13}\\
0 & 0 &g_{23}\\
0&0&0\end{pmatrix}.
\]
Since $\fv$ is a Lie subalgebra, it also contains
\[
[f,g]=f_{12}\circ g_{23},
\]
where $f_{12}\circ g_{23}$ is considered as an element of $W_{-1}\End_\mathbb{K}(\omega M)$ via the natural embedding of $\Hom_\mathbb{K}(\omega A_3,\omega A_1)$ (as the (13)-entry). The claim now follows because the map
\[
\Hom_\mathbb{K}(\omega A_2,\omega A_1)\otimes \Hom_\mathbb{K}(\omega A_3,\omega A_2) \rightarrow \Hom_\mathbb{K}(\omega A_3,\omega A_1)
\]
defined by $f_{12}\otimes g_{23}\mapsto f_{12}\circ g_{23}$ is surjective.

In view of \eqref{eq10}, we thus obtain:
\begin{prop}\label{prop: u maximal iff u-1 maximal}
Suppose $A_2\neq 0$. We have $\fu(M)=W_{-1}\inEnd(M)$ if and only if 
\[
\fu_{-1}(M)=\inHom(A_2,A_1)\oplus \inHom(A_3,A_2).
\]
\end{prop}

\begin{rem}
This is a generalization of Corollary 4.6 of \cite{BP1}.
\end{rem}

\begin{rem}
Most of the constructions of \S \ref{sec: initial considerations} can be generalized without much difficulty to filtrations with an arbitrary number of steps in a tannakian category over a field of characteristic zero.\footnote{I thank an anonymous referee for pointing this out.} We focused on the case of blended extensions because the remainder of the results of the paper are for this case.
\end{rem}

\section{The determination of $\fu_{-1}(M)$}\label{sec: char of u-1}
Throughout this section we assume that our data of \S \ref{sec: data} satisfies condition (C1) of \S \ref{sec: conditions C1 and C2}.

\subsection{}
A key tool for the proof of parts (b),(f), and (g) of Theorem \ref{thm: intro} is the following proposition.

\begin{prop}\label{prop: map for Ext groups}
(a) For every semisimple object $X$ of $\langle M\rangle^{\otimes}$ (where $M$ satisfies condition (C1) of \S \ref{sec: conditions C1 and C2}), there is a map
\begin{equation}\label{eq43}
\Ext^1_{\langle M\rangle^{\otimes}}(\mathbbm{1}, X) \xrightarrow{\ \simeq \ } \Hom(\fu^{\ab}(M), X)
\end{equation}
such that the image of the class of an extension $E$ of $\mathbbm{1}$ by $X$ in $\langle M\rangle^{\otimes}$ under \eqref{eq43} is calculated as follows: Choose a fiber functor $\omega$ and a linear section of $\omega E\twoheadrightarrow \omega\mathbbm{1}$ to decompose $\omega E \cong \omega X\oplus \omega \mathbbm{1}$ as vector spaces, express the action of $\mathcal{U}(M,\omega)$ on $\omega E$ with respect to this decomposition to obtain a morphism of algebraic groups $\mathcal{U}(M,\omega)\rightarrow \omega X$ over $\mathbb{K}$ (where $\omega X$ is a vector group), and then take the logarithm of this morphism. The map $\omega \fu^{\ab}(M)\rightarrow \omega X$ obtained after passing to the abelianization of $\omega \fu(M)$ is the image under $\omega$ of the morphism $\fu^{\ab}(M)\rightarrow X$ corresponding to the class of $E$ under \eqref{eq43}.

(b) The map \eqref{eq43} is a $\mathbb{K}$-linear isomorphism that is functorial with respect to morphisms $X\rightarrow X'$ between semisimple objects $X,X'$ of $\langle M\rangle^{\otimes}$. Moreover, it is independent of the choice of $\omega$.
\end{prop}

\begin{proof}
A more general statement is proven in \cite[\S 2]{Es26} (see also therein for the relationship to the inflation-restriction sequence in group cohomology). Here, for the sake of completeness, we include an argument for the statement claimed above. Also, the approach taken below to construct the map \eqref{eq43} makes the naturality of the map clearer than the approach in \cite{Es26}.

Since $X$ is semisimple, it belongs to $\langle A_1,A_2,A_3\rangle^{\otimes}$. Starting with an extension
\begin{equation}\label{eq12}
\begin{tikzcd}
   0 \arrow[r] & X \arrow[r, ] & E \arrow[r, ] & \mathbbm{1}  \arrow[r] & 0
\end{tikzcd}
\end{equation}
in $\langle M\rangle^{\otimes}$, the inclusion of $\langle E\rangle^{\otimes}$ in $\langle M\rangle^{\otimes}$ gives a morphism $\fg(M)\rightarrow \fg(E)$ of Lie algebra objects in $\bT$. Since $X$ belongs to $\langle A_1,A_2,A_3\rangle^{\otimes}$, the morphism $\fg(M)\rightarrow \fg(E)$ restricts to a morphism of Lie algebra objects $\fu(M)=\fg(M,A_1\oplus A_2\oplus A_3)\rightarrow \fg(E,X)$. This morphism only depends on the object $E$, rather than the extension \eqref{eq12}. The choice of the extension will give an embedding $\fg(E,X)\hookrightarrow \inHom(\mathbbm{1},X)\cong X$ of Lie algebra objects, where $\inHom(\mathbbm{1},X)$ is considered as an abelian Lie subobject of $\inEnd(E)$ (see for instance, \cite{EM1}). This morphism is easily seen to only depend on the class of \eqref{eq12} in $\Ext^1(\mathbbm{1},X)$. Since $X$ is abelian, the Lie algebra morphism $\fu(M)\rightarrow \fg(E,X) \rightarrow X$ factors through a morphism $\fu(M)^{\ab}\rightarrow X$. The map \eqref{eq43} sends the class of \eqref{eq12} to this morphism $\fu(M)^{\ab}\rightarrow X$.\footnote{I thank an anonymous referee for suggesting this outline for the construction of the map \eqref{eq43}, which makes the naturality of the construction clearer.} That the latter morphism is obtained by the recipe outlined in the statement of part (a) follows from the definitions of the two maps $\fu(M)\rightarrow \fg(E,X)$ and $\fg(E,X)\hookrightarrow X$.

The functoriality of \eqref{eq43} in $X$ can be checked using the construction of the map and pushforwards of extensions. The linearity follows from functoriality.

Finally, one way to see that \eqref{eq43} is bijective is to construct its inverse as follows. Fix a fiber functor $\omega$ and a Levi factor $\mathcal{R}$ of $\mathcal{G}(M,\omega)$. Consider the exact sequence of $\mathcal{R}$-representations
\begin{equation}\label{eq14}
\begin{tikzcd}
   0 \arrow[r] & \omega X \arrow[r, ] & \omega X\oplus \omega\mathbbm{1}  \arrow[r, ] & \omega\mathbbm{1}  \arrow[r] & 0
\end{tikzcd}
\end{equation}
with the direct sum action and the obvious inclusion and projection maps. Let $f: \fu^{\ab}(M)\rightarrow X$ be a morphism in $\mathbf{T}$. Then the composition $\omega\fu(M)\twoheadrightarrow \omega\fu(M)^{\ab}\xrightarrow{\omega f} \omega X$, which we also denote by $\omega f$, is a $\mathcal{G}(M,\omega)$-equivariant morphism of Lie algebras. It lifts to a $\mathcal{G}(M,\omega)$-equivariant morphism of algebraic groups $\exp(\omega f): \mathcal{U}(M,\omega)\rightarrow \omega X$ (the action of $\mathcal{G}(M,\omega)$ on $\mathcal{U}(M,\omega)$ being by conjugation). We use $\exp(\omega f)$ in the natural way to define a $\mathcal{U}(M,\omega)$-action on $\omega X\oplus \omega\mathbbm{1}$: an element $u\in \mathcal{U}(M,\omega)$ will act on $\omega X\oplus \omega\mathbbm{1}$ by left multiplication by the matrix
\[
\begin{pmatrix}
\Id_{\omega X} & \exp(\omega f)(u)\\
0 & \Id_{\omega\mathbbm{1}}\end{pmatrix}\in \GL(\omega X\oplus \omega \mathbbm{1}).
\]
This makes \eqref{eq14} a sequence of $\mathcal{U}(M,\omega)$-representations. The $\mathcal{R}$-equivariance of $\exp(\omega f)$ guarantees that the action of $\mathcal{U}(M,\omega)$ on $\omega X\oplus \omega\mathbbm{1}$ defined above glues with the direct sum action of $\mathcal{R}$ to define a $\mathcal{G}(M,\omega)$-action on $\omega X\oplus \omega\mathbbm{1}$. The sequence \eqref{eq14} becomes an extension of $\mathcal{G}(M,\omega)$-representations, and hence isomorphic to the class of $\omega E$ for some extension $E$ of $\mathbbm{1}$ by $X$ in $\langle M\rangle^{\otimes}$. One can see using reductivity of $\mathcal{R}$ that sending $f$ to the class of $E$ in $\Ext^1_{\langle M\rangle^{\otimes}}(\mathbbm{1},X)$ we obtain the inverse of \eqref{eq43}.
\end{proof}

\subsection{}\label{sec: injection Ext(1,X)->Hom(u,X)} Let $X$ be an object of $\langle M\rangle^{\otimes}$. The map $\fu(M)\twoheadrightarrow \fu^{\ab}(M)$ induces an injection
\[
\Hom(\fu^{\ab}(M), X) \hookrightarrow \Hom(\fu(M), X)
\]
that is functorial in $X$. If $X$ is semisimple, composing this with the isomorphism \eqref{eq43} we obtain an injection
\begin{equation}\label{eq18}
\Psi_X : \ \Ext^1_{\langle M\rangle^{\otimes}}(\mathbbm{1}, X) \hookrightarrow \Hom(\fu(M), X)
\end{equation}
that is also functorial in $X$. The image of an extension $E$ under this can be calculated just like the image under \eqref{eq43} of $E$ following the procedure outlined in Proposition \ref{prop: map for Ext groups}(a), except we skip the final step of passing to the map induced on the abelianization of $\omega \fu(M)$.

\subsection{}\label{subsection: proof of part b of the main thm}
Here we establish the characterization of $\fu_{-1}(M)$ given in Theorem \ref{thm: intro}(b). We shall continue to assume that $A_1,A_2,A_3$ are semisimple. Following \S \ref{sec:introduction}, we shall set
\[
V:=\inHom(A_2,A_1)\oplus \inHom(A_3,A_2).
\]
Recall from \S \ref{sec: data} that $\mathcal{L}$ and $\mathcal{N}$ denote, respectively, the elements of $\Ext^1(\mathbbm{1},\inHom(A_2,A_1))$ and $\Ext^1(\mathbbm{1},\inHom(A_3,A_2))$) corresponding to $L$ and $N$. Both $\mathcal{L}$ and $\mathcal{N}$ belong to the subgroups $\Ext^1_{\langle M\rangle^{\otimes}}$ of the corresponding $\Ext^1$ groups. Recall that we wrote the coordinates of the canonical map $\pi$ of \S \ref{sec: construction of pi} as $\pi_{12}$ and $\pi_{23}$. Denote the composition
\[
\fu(M)\hookrightarrow
\begin{tikzcd}[column sep = .8in]
W_{-1}\inEnd(M) \arrow[r, twoheadrightarrow, "\text{$\pi=(\pi_{12},\pi_{23})$}"] & V
\end{tikzcd}
\]
also by $\pi=(\pi_{12},\pi_{23})$. We will denote the injection $\Psi_X$ of \eqref{eq18} for the choices $X=\inHom(A_2,A_1)$ and $X=\inHom(A_3,A_2)$ respectively by $\Psi_{12}$ and $\Psi_{23}$.

The main ingredient of the proof of Theorem \ref{thm: intro}(b) is the following proposition.

\begin{prop}\label{prop: calculation of image of L and N under the canonical isomorphism Psi}
We have
\begin{align*}
\Psi_{12}(\sL) & = \pi_{12}\\
\Psi_{23}(\sN) & = \pi_{23}.
\end{align*}
\end{prop}
\begin{proof}
The two formulas are verified similarly, so we will only prove the second one. We will first explicitly write a representative of $\sN$, and then follow the procedure outlined in Proposition \ref{prop: map for Ext groups}(a) to calculate $\Psi_{23}(\sN)$. 

The extension class $\mathcal{N}$ is represented by the extension
\[
\begin{tikzcd}
0 \arrow[r] & \inHom(A_3,A_2) \arrow[r] & \inHom(A_3,N)^\dagger \arrow[r] & \mathbbm{1} \arrow[r] & 0,
\end{tikzcd}
\]
also denoted by $\sN$, which after applying a fiber functor $\omega$ becomes the extension
\[
\begin{tikzcd}
0 \arrow[r] & \Hom_\mathbb{K}(\omega A_3,\omega A_2) \arrow[r] & \Hom_\mathbb{K}(\omega A_3,\omega N)^\dagger \arrow[r] & \mathbb{K} \arrow[r] & 0.
\end{tikzcd}
\]
Here, $\inHom(A_3,N)^\dagger$ is the subobject of $\inHom(A_3,N)$ whose image
\[
\Hom_\mathbb{K}(\omega A_3,\omega N)^\dagger \subset \Hom_\mathbb{K}(\omega A_3,\omega N)
\]
under $\omega$ consists of all linear maps $g: \omega A_3\rightarrow \omega N$ such that the composition $\omega A_3\xrightarrow{g} \omega N \twoheadrightarrow \omega A_3$ (the latter map being the structure map) is a scalar multiple, denoted by $\lambda(g)\in \mathbb{K}$, of the identity map on $\omega A_3$. The surjective arrow in $\omega \sN$ is the map $\lambda: g\mapsto \lambda(g)$, and the injective arrow in $\omega \sN$ is given by the inclusion $\omega A_2\hookrightarrow \omega N$. See \cite[\S 3.2]{EM2} for more details.

Choose a splitting of $\omega M$ in the sense of \S \ref{sec: splitting}. This splitting induces a section for the exact sequence $\omega \sN$ above. We identify $\Hom_\mathbb{K}(\omega A_3,\omega N)^\dagger$ with $\Hom_\mathbb{K}(\omega A_3,\omega A_2)\oplus \mathbb{K}$ via this section. 

Let $f=(f_{ij})\in \omega\fu(M)\subset \End_\mathbb{K}(\omega M)$, where we write endomorphisms of $\omega M$ as $3\times 3$ matrices using our splitting of $\omega M$. Set $\sigma:=\exp(f)=1+f+f^2/2$. Note that $\sigma$ has the same super-diagonal entries as $f$ (i.e., $f_{12}$ and $f_{23}$). The element $\sigma\in \mathcal{U}(M,\omega)$ acts on elements of $\omega \inHom(A_3,N)$ by sending $g: \omega A_3\rightarrow \omega N$ to $\sigma_N g\sigma_{A_3}^{-1}$, where $\sigma_N$ and $\sigma_{A_3}$ are the actions of $\sigma$ on $\omega N$ and $\omega A_3$, respectively. Writing elements of $\GL(\omega N)$ as $2\times 2$ matrices via the isomorphism $\omega N\cong \omega A_2\oplus \omega A_3$ given by our splitting, given any $g=(g_{23},\lambda)\in \Hom_\mathbb{K}(\omega A_3,\omega N)^\dagger$, since $\sigma_{A_3}=1$ and $\sigma_N$ is the bottom right $2\times 2$ submatrix of $\sigma$, we have
\[
\sigma\cdot g  = \sigma_N g \sigma_{A_3}^{-1} = \begin{pmatrix} 1& f_{23}\\ 0& 1\end{pmatrix}
\begin{pmatrix} g_{23} \\ \lambda\end{pmatrix}.
\]
Thus the image of $\Psi_{23}(\sN): \fu(M)\rightarrow \inHom(A_3,A_2)$ under $\omega$ is the logarithm of the morphism 
\[
\mathcal{U}(M,\omega)\rightarrow \Hom_\mathbb{K}(\omega A_3,\omega A_2) \hspace{.3in} \sigma=\exp(f) \mapsto f_{23} 
\] 
of algebraic groups. Thus $\Psi_{23}(\sN)$ after applying $\omega$ is simply the projection $f\mapsto f_{23}$, as desired.
\end{proof}

We can now deduce our characterization of $\fu_{-1}(M)$.

\begin{thm}\label{thm: characterization of u-1}
Suppose that $M$ satisfies condition (C1) of \S \ref{sec: conditions C1 and C2}. Consider $(\sL, \sN)$ as an element of $\Ext^1(\mathbbm{1}, V)$ via the canonical isomorphism
\[
\Ext^1(\mathbbm{1}, \inHom(A_2,A_1))\oplus \Ext^1(\mathbbm{1}, \inHom(A_3,A_2)) \cong \Ext^1(\mathbbm{1}, V).
\]
Then $\fu_{-1}(M)$ is the smallest subobject of $V$ with the property that the image of $(\sL,\sN)$ in $\Ext^1(\mathbbm{1}, V/\fu_{-1}(M))$ is split.
\end{thm}

\begin{proof}
By the functoriality property of the map $\Psi_X$ of \eqref{eq18}, we have
\[
\Psi_V(\sL, \sN)=(\Psi_{12}(\sL),\Psi_{23}(\sN)\stackrel{\text{Prop. \ref{prop: calculation of image of L and N under the canonical isomorphism Psi}}}{=} \pi.
\]
Let $X\subset V$. Again by functoriality, we have a commutative diagram
\[
\begin{tikzcd}[column sep = large]
\Ext^1_{\langle M\rangle^{\otimes}}(\mathbbm{1}, V) \ar[d, "q_\ast"] \arrow[r, hookrightarrow, "\Psi_V"] & \Hom(\fu(M), V) \ar[d, "q\circ -"] \\
\Ext^1_{\langle M\rangle^{\otimes}}(\mathbbm{1}, V/X) \arrow[r, hookrightarrow, "\Psi_{V/X}"] & \Hom(\fu(M), V/X),
\end{tikzcd}
\]
where the vertical maps are given by the quotient map $q: V\rightarrow V/X$. The object $\fu_{-1}(M)=\pi(\fu(M))$ is contained in $X$ if and only if $q\circ\pi$ vanishes on $\fu(M)$, which in view of the commutativity of the diagram is equivalent to $\Psi_{V/X}\circ q_\ast(\sL, \sN)=0$. We obtain the result since $\Psi_{V/X}$ is injective.
\end{proof}

\begin{rem}
Since $V$ is assumed to be semisimple, every subobject of $V$ is the kernel of an endomorphism of $V$. By a similar argument to above, the subobject $\fu_{-1}(M)$ of $V$ is the intersection of the kernels of all the endomorphisms $\phi$ of $V$ such that $\phi_\ast(\sL,\sN) = 0$.
\end{rem}

\begin{rem}
The reader can refer to \cite[Theorem 1.1.1]{Es26} for a generalization of Theorem \ref{thm: characterization of u-1} to filtrations with an arbitrary number of steps and possibly non-semisimple graded pieces, as well as some applications (particularly, regarding maximality of $\fu(M)$).
\end{rem}

\section{The characterizations of $\fu_{-2}(M)$ and its subquotients}\label{sec: characterization of u-2}
Our first goal in this section is to give a characterization of $\fu_{-2}(M)$ in the spirit of the characterization of the tannakian group of an extension given in \cite[Theorem 3.3.1]{EM1}. This will be the subject of \S \ref{sec: lem for charaterization of u-2} - \S \ref{sec: general characterization of u-2} below, and is in the full generality of the data of \S \ref{sec: data}, without necessarily satisfying any of the conditions of \S \ref{sec: conditions C1 and C2}. After that, we will turn our focus to the derived Lie algebra $[\fu(M),\fu(M)]\subset \fu_{-2}(M)$ and the quotient $\fu_{-2}(M)/ [\fu(M),\fu(M)]$.

\subsection{}\label{sec: lem for charaterization of u-2}
Recall from \S \ref{sec: define g(X,Y)} that for every objects $X$ of $\bT$ and $Y$ of $\langle X\rangle^{\otimes}$ the notation $\fg(X,Y)$ means the canonical subobject of $\inEnd(X)$ (or $\fg(X)$) whose image under any fiber functor $\omega$ is the Lie algebra of the kernel of the natural surjection $\mathcal{G}(X,\omega)\rightarrow \mathcal{G}(Y,\omega)$. We denoted this kernel by $\mathcal{G}(X,Y,\omega)$.

\begin{lemma}\label{lem: u-2=lie(ker(G(M)-> G(L+N)))}
We have
\[
\fu_{-2}(M) = \fg(M, L\oplus N).
\]
\end{lemma} 
\begin{proof}
Let $\omega$ be a fiber functor. Since $A_1\oplus A_2\oplus A_3$ belongs to $\langle L\oplus N\rangle^{\otimes}$, we have
\[
\mathcal{G}(M,L\oplus N,\omega) \subset \mathcal{G}(M,A_1\oplus A_2\oplus A_3,\omega).
\]
In particular, $\mathcal{G}(M,L\oplus N,\omega)$ is unipotent and $\omega \fg(M, L\oplus N)\subset \omega \fu(M)$.

We show that $\omega \fu_{-2}(M)$ and $\omega \fg(M, L\oplus N)$ coincide in $W_{-1}\End_\mathbb{K}(\omega M)$. Choose a splitting of $\omega M$ (see \S \ref{sec: splitting}) to express elements of $\GL(\omega M)$ and $\End_\mathbb{K}(\omega M)$ as $3\times 3$ matrices. Let $f=(f_{ij})\in \omega\fu_{-2}(M)$ (so $f_{13}$ is the corresponding map $\omega A_3\rightarrow \omega A_1$ and the rest of the $f_{ij}$ are all zero). Then
\begin{equation}\label{eq6}
\exp(f)=\begin{pmatrix}1&0&f_{13}\\0&1&0\\0&0&1\end{pmatrix}\in \mathcal{G}(M,\omega)
\end{equation}
acts trivially on both $L$ and $N$. Thus $\omega \fu_{-2}(M)\subset \omega \fg(M, L\oplus N)$.

On the other hand, if $f\in \omega \fg(M, L\oplus N)$, then $\exp(f)\in \mathcal{G}(M, L\oplus N, \omega)$ acts trivially on $\omega L$ and $\omega N$ and hence is of the form \eqref{eq6} for some $f_{13}$. That is, $f$ belongs to the subspace $\Hom_\mathbb{K}(\omega A_3,\omega A_1)\hookrightarrow W_{-1}\End_\mathbb{K}(\omega M)$.
\end{proof}

\subsection{}
Let $X$ be an object of $\bT$ and $\mathcal{E}\in \Ext^1(\mathbbm{1},X)$. In \cite{EM2} with Murty we defined what it means for $\mathcal{E}$ to originate from a full tannakian subcategory $\bS$ of $\bT$ that is closed under taking subobjects (or subquotients): we say $\mathcal{E}$ originates from $\bS$ if there exist a subobject $X'$ of $X$ in $\bS$ and an extension $\mathcal{E'}\in \Ext^1_\bS(\mathbbm{1},X')$ that pushes forward to $\sE$. 

In what follows, we shall only deal with instances of the above definition that $X$ is in $\bS$. In that case, the definition simplifies: $\sE$ originates from $\bS$ if and only if $\sE$ is in the subgroup
\[
\Ext^1_\bS(\mathbbm{1},X) \subset \Ext^1(\mathbbm{1},X)
\]
if and only if the middle object of $\sE$ belongs to $\bS$.

\subsection{}\label{sec: general characterization of u-2}
Let
\[\sM^{h}\in \Ext^1(\mathbbm{1},\inHom(N,A_1))\]
be the element corresponding to the horizontal extension on the middle row of \eqref{eq1} under the canonical isomorphism
\[
\Ext^1(N,A_1) \cong \Ext^1(\mathbbm{1},\inHom(N,A_1)).
\]
We are ready to give our characterization of $\fu_{-2}(M)$.

\begin{thm}\label{thm: characterization of u-2}
Let $\fv$ be a subobject of $\inHom(A_3,A_1)$. Then the extension 
\[\mathcal{M}^h/\fv \, \in \, \Ext^1(\mathbbm{1}, \inHom(N,A_1)/\fv)\]
originates from $\langle L\oplus N \rangle^\otimes$ if and only if $\fu_{-2}(M)\subset \fv$. That is, $\fu_{-2}(M)$ is the smallest subobject of $\inHom(A_3,A_1)$ such that $\mathcal{M}^h/\fu_{-2}(M)$ originates from $\langle L\oplus N \rangle^\otimes$.
\end{thm}

\begin{proof}
We will prove that for every subobject $\fv$ of $\inHom(N,A_1)$, the extension $\mathcal{M}^h/\fv$ originates from $\langle L\oplus N \rangle^\otimes$ if and only if $\fu_{-2}(M)\subset \fv$. We may assume that $N\neq 0$. Fix a fiber functor $\omega$. The extension $\sM^h$ is given by 
\[
\begin{tikzcd}
0 \arrow[r] & \inHom(N,A_1) \arrow[r] & \inHom(N,M)^\dagger \arrow[r] & \mathbbm{1} \arrow[r] & 0,
\end{tikzcd}
\]
where the notation is in line with our notation in the proof of Proposition \ref{prop: calculation of image of L and N under the canonical isomorphism Psi}. That is, the image of $\inHom(N,M)^\dagger$ under $\omega$ consists of linear maps $g:\omega N\rightarrow \omega M$ whose composition with the structure map $\omega M\twoheadrightarrow \omega N$ is a scalar multiple $\lambda(g)$ of the identity map on $\omega N$. After applying $\omega$, the map $\inHom(N,M)^\dagger\rightarrow \mathbbm{1}$ sends $g$ to $\lambda(g)$.

By \cite[Lemma 3.4.2]{EM2}, $\sM^h/\fv$ originates from $\langle L\oplus N\rangle^{\otimes}$ if and only if $\omega(\mathcal{M}^h/\fv)$ splits in the category of representations of $\mathcal{G}(M,L\oplus N,\omega)$. Since $\inHom(N,A_1)/\fv$ belongs to $\langle L\oplus N\rangle^{\otimes}$, the action of $\mathcal{G}(M,L\oplus N,\omega)$ on $\omega(\inHom(N,A_1)/\fv)$ is trivial. It follows that $\mathcal{M}^h/\fv$ originates from $\langle L\oplus N\rangle^{\otimes}$ if and only if, identifying $\inHom(N,M)^\dagger\subset W_{-1}\inEnd(M)$ and choosing a splitting of $\omega M$ to identify $\omega M=\omega A_1\oplus \omega A_2\oplus \omega A_3$ as vector spaces and hence writing maps $\omega M\rightarrow \omega M$ as matrices, the element
\begin{equation}\label{eq8}
\begin{pmatrix}
0 & & \\
& 1 & \\
&& 1\end{pmatrix} +\omega \fv \ \in \ \Hom_\mathbb{K}(\omega N,\omega M)^\dagger/\omega \fv \ \subset \ W_{-1}\End_\mathbb{K}(\omega M)/\omega \fv
\end{equation}
is fixed by $\mathcal{G}(M,L\oplus N,\omega)$. By Lemma \ref{lem: u-2=lie(ker(G(M)-> G(L+N)))}  and the fact that $\mathcal{G}(M,L\oplus N,\omega)$ is unipotent, 
\[
\mathcal{G}(M,L\oplus N,\omega) = \exp(\omega \fu_{-2}(M)).
\]
For every $f\in \Hom_\mathbb{K}(\omega A_3,\omega A_1)\subset W_{-1}\End_\mathbb{K}(\omega M)$, a direct computation gives
\[
\exp(f)\begin{pmatrix}
0 & & \\
& 1 & \\
&& 1\end{pmatrix}\exp(f)^{-1} - \begin{pmatrix}
0 & & \\
& 1 & \\
&& 1\end{pmatrix} = f.
\]
It follows that $\mathcal{G}(M,L\oplus N,\omega)$ fixes \eqref{eq8} if and only if $\omega \fu_{-2}(M)\subset \omega \fv$. 
\end{proof}

Thus roughly speaking, $\fu_{-2}(M)$ captures the obstruction to the extension $\sM^{h}$ being originated from $\langle L\oplus N\rangle^{\otimes}$, or rather, the obstruction to $M$ being an object of $\langle L\oplus N\rangle^{\otimes}$. In particular, one has the following corollary:

\begin{cor}\label{cor: u-2 is zero iff M is generated by L and N}
We have $\fu_{-2}(M)=0$ if and only if $M$ is contained in $\langle L\oplus N\rangle^{\otimes}$.
\end{cor}
\begin{proof}
By Theorem \ref{thm: characterization of u-2}, the object $\fu_{-2}(M)$ is zero if and only if $\inHom(N,M)^\dagger$ (i.e. the middle object of $\sM^h$, see the proof of Theorem \ref{thm: characterization of u-2}) belongs to the subcategory $\langle L\oplus N\rangle^{\otimes}$. The claim now follows on noting that $\langle L\oplus N\rangle^{\otimes}$ contains $M$ if and only if it contains $\inHom(N,M)^\dagger$. (Note for the {\it if} implication: For nonzero $N$ the evaluation morphism $\inHom(N,M)^\dagger\otimes N\rightarrow M$ is surjective.)
\end{proof}

\begin{rem}
One can give an analogous characterization of $\fu_{-2}(M)$ in terms of the extension
\[\sM^{v}\in \Ext^1(\mathbbm{1},\inHom(A_3,L))\]
corresponding to the middle column of \eqref{eq1}. By a similar argument to the one above, $\fu_{-2}(M)$ is the smallest subobject of $\inHom(A_3,A_1)$ such that the pushforward 
\[\sM^v/\fu_{-2}(M)\in \Ext^1(\mathbbm{1},\inHom(A_3,L)/\fu_{-2}(M))\]
originates from $\langle L\oplus N\rangle^{\otimes}$.
\end{rem}

\subsection{}\label{sec: explicit characterization of [u,u]}
Our subject of study in the next two subsections is the derived algebra $[\fu(M),\fu(M)]$. In this subsection we shall see that this subobject of $\inHom(A_3,A_1)$ is completely determined by $\fu_{-1}(M)$; in fact, we will see that one can very easily compute $[\fu(M),\fu(M)]$ from $\fu_{-1}(M)$. In particular, when condition (C1) of \S \ref{sec: conditions C1 and C2} holds, combining with Theorem \ref{thm: intro}(b) we will have an explicit description of the derived algebra. We do not need to assume any conditions for the discussion of this subsection.

As before, let $V$ be $\inHom(A_2,A_1)\oplus \inHom(A_3,A_2)$. Let
\begin{equation}\label{eq21}
V\otimes V \rightarrow \inHom(A_3,A_1)
\end{equation}
be the morphism that after applying a fiber functor $\omega$ is given by
\[
(f_{12},f_{23})\otimes (g_{12},g_{23}) \mapsto f_{12}\circ g_{23}-g_{12}\circ f_{23}.
\]
That this is a morphism is easily seen by verifying compatibility with the actions of $\mathcal{G}(A_1\oplus A_2\oplus A_3, \omega)$. We use the notation $\{ \, ,\}$ for the map \eqref{eq21} and the associated pairing.

One has a commutative diagram
\[
\begin{tikzcd}[column sep = small]
\fu(M)\otimes \fu(M) \arrow[rd, "\text{$[\, , ]$}"] \arrow[rr, twoheadrightarrow, "\pi\otimes\pi "] & & \fu_{-1}(M)\otimes \fu_{-1}(M) \arrow[ld, "\text{$\{ \, , \}$}"] \\
& \inHom(A_3,A_1). &
\end{tikzcd}
\]
Indeed, choose a fiber functor $\omega$ and a splitting of $\omega M$ to identify $\omega M$ as $\omega A_1\oplus \omega A_2\oplus \omega A_3$. Recall that $\omega \pi$ simply sends an element $f=(f_{ij})$ of $\fu(M)\subset W_{-1}\End_\mathbb{K}(\omega M)$ to $(f_{12},f_{23})$. The commutativity of the diagram above after applying $\omega$ is checked by a direct computation.

We thus obtain the following:
\begin{prop}\label{prop: [u,u]={u-1,u-1}}
The derived algebra $[\fu(M),\fu(M)]$ is the image of $\fu_{-1}(M)\otimes \fu_{-1}(M)$ under the map \eqref{eq21}.
\end{prop}
Combining with \S \ref{sec: u-1 is in u(L)+u(N)} we obtain the following corollary:
\begin{cor}
The subobject $[\fu(M),\fu(M)]$ of $\inHom(A_3,A_1)$ only depends on the extensions $L$ and $N$ (and not on the particular blended extension $M$). In particular, whether or not the Lie algebra $\fu(M)$ of a blended extension $M$ of $N$ by $L$ is abelian only depends on the extensions $N$ and $L$.
\end{cor}

\subsection{}
The last subsection allows us to explicitly calculate the derived algebra $[\fu(M),\fu(M)]$ from $\fu_{-1}(M)$. Here we record a more conceptual (although possibly less practical) characterization of $[\fu(M),\fu(M)]$ under condition (C1). The characterization is rather basic and surely well known; it is included here mainly for the purpose of completeness of the discussion.

\begin{prop}\label{prop: [u,u]}
Assume condition (C1) of \S \ref{sec: conditions C1 and C2}. For every subobject $\fv$ of $\fu_{-2}(M)$ the following statements are equivalent:
\begin{itemize}
\item[(i)] The pushforward of the extension
\begin{equation}\label{eq19}
\begin{tikzcd}
0 \arrow[r] & \fu_{-2}(M) \arrow[r] & \fu(M) \arrow[r] & \fu_{-1}(M) \arrow[r] & 0
\end{tikzcd}
\end{equation}
along the quotient map $\fu_{-2}(M)\twoheadrightarrow \fu_{-2}(M)/\fv$ splits.
\item[(ii)] The quotient $\fu(M)/\fv$ is a semisimple object of $\bT$.
\item[(iii)] We have $[\fu(M),\fu(M)] \subset \fv$.
\end{itemize}
In particular, $[\fu(M),\fu(M)]$ is zero if and only if $\fu(M)$ is a semisimple object of $\bT$.
\end{prop}

\begin{proof}
The equivalence of (i) and (ii) is clear upon recalling that thanks to condition (C1), $\fu_{-1}(M)$ and $\fu_{-2}(M)$ are semisimple. We will argue that (ii) and (iii) are equivalent. The object $\fu(M)/\fv$ is semisimple if and only if the adjoint action of $\mathcal{G}(M,\omega)$ on it (after applying a fiber functor $\omega$) factors through an action of $\mathcal{G}(A_1\oplus A_2\oplus A_3, \omega)$, or equivalently, if and only if the action of $\mathcal{G}(M, A_1\oplus A_2\oplus A_3, \omega)$ on it is trivial. The latter statement is equivalent to the triviality of the action of $\omega \fu(M)$ on $\omega \fu(M)/\omega \fv$ induced by the Lie bracket, which is in turn equivalent to the statement that $\omega \fv$ contains $[\omega \fu(M),\omega \fu(M)]$.

The final assertion follows from specializing to the case $\fv=0$.
\end{proof}
Thus $[\fu(M),\fu(M)]$ is the smallest subobject of $\fu_{-2}(M)$ such that $\fu(M)/[\fu(M),\fu(M)]$ is a semisimple object of $\bT$, or equivalently, such that the pushforward of the extension of \eqref{eq19} along $\fu_{-2}(M)\twoheadrightarrow \fu_{-2}(M)/[\fu(M),\fu(M)]$ splits. Since every subobject of $\fu(M)$ is a Lie subobject (in fact, a Lie ideal subobject), we may replace the term {\it the smallest subobject} in the previous sentence by {\it the smallest Lie subobject}.

\subsection{}\label{sec: char of u-2/[u,u]}

We now turn our attention to the quotient of $\fu_{-2}(M)$ by the derived algebra. The following result can be useful for characterizing this quotient when $M$ satisfies conditions (C1) and (C2) of \S \ref{sec: conditions C1 and C2}.

\begin{prop}\label{prop: characterization of u-2/[u,u]}
Assume (C1) and (C2). There is a canonical isomorphism
\begin{equation}\label{eq20}
\Ext^1_{\langle M\rangle^{\otimes}}(A_3,A_1) \ \cong \ \Hom(\frac{\fu_{-2}(M)}{[\fu(M),\fu(M)]}, \inHom(A_3,A_1)).
\end{equation}
\end{prop}

\begin{proof}
Thanks to (C1) and Proposition \ref{prop: map for Ext groups} (applied with $X=\inHom(A_3,A_1)$), we have a canonical isomorphism
\[
\Ext^1_{\langle M\rangle^{\otimes}}(A_3,A_1) \cong \ \Ext^1_{\langle M\rangle^{\otimes}}(\mathbbm{1},\inHom(A_3,A_1)) \cong \Hom(\fu^{\ab}(M), \inHom(A_3,A_1)).
\]
Since $\fu^{\ab}(M)$ is a semisimple object of $\bT$, the restriction map 
\[
\Hom(\fu^{\ab}(M), \inHom(A_3,A_1)) \rightarrow \Hom(\frac{\fu_{-2}(M)}{[\fu(M),\fu(M)]}, \inHom(A_3,A_1))
\]
induced by the inclusion of $\fu_{-2}(M)/[\fu(M),\fu(M)]$ in $\fu^{\ab}(M)$ is surjective. Condition (C2) guarantees that this map is an isomorphism, as there are no nonzero morphisms from $\fu_{-1}(M)$ to $\inHom(A_3,A_1)$. Composing this isomorphism with the previous isomorphism we obtain \eqref{eq20}.
\end{proof}

The object $\fu_{-2}(M)/[\fu(M),\fu(M)]$ is a subquotient of $\inHom(A_3,A_1)$, and hence assuming (C1), also a subobject of it. Thus the last proposition has the following immediate consequence:
\begin{cor}\label{cor: when is [u,u]=u-2}
Assume (C1) and (C2). Then $[\fu(M),\fu(M)]=\fu_{-2}(M)$ if and only if 
\[\Ext^1_{\langle M\rangle^{\otimes}}(A_3,A_1)=0.\]
In particular, one has $\fu_{-2}(M)=0$ if and only if $\fu(M)$ is an abelian Lie algebra and the group $\Ext^1_{\langle M\rangle^{\otimes}}(A_3,A_1)$ vanishes.
\end{cor}

\subsection{}\label{sec: at most one M with vanishing u-2}
Building on the last corollary, we end this section with another consequence of Proposition \ref{prop: characterization of u-2/[u,u]}, which is an observation on blended extensions with vanishing $\fu_{-2}$.
\begin{prop}\label{prop: at most one M with vanishing u-2}
Assume conditions (C1) and (C2) of \S \ref{sec: conditions C1 and C2}. Fixing $L$ and $N$, there exists at most one blended extension $M$ in $\Extpan(N,L)$ such that $\fu_{-2}(M)$ vanishes.
\end{prop}

\begin{proof}
Suppose that there exist blended extensions $M_1$ and $M_2$ of $N$ by $L$ both with vanishing $\fu_{-2}$. Then in particular, $\fu(M_1)$ and $\fu(M_2)$ are abelian and hence as objects of $\bT$, they are semisimple (see Proposition \ref{prop: [u,u]}). There is a natural injection
\[
\fu(M_1\oplus M_2) \hookrightarrow \fu(M_1)\oplus \fu(M_2),
\]
induced by the natural embedding of the tannakian group of $M_1\oplus M_2$ in the product of the tannakian groups of $M_1$ and $M_2$. Here, $\fu(M_1\oplus M_2)$ is the object $\fu$ of the direct sum blended extension $M_1\oplus M_2$ (a blended extension of $N^2$ by $L^2$), defined according to \S \ref{sec: definition of u}; thanks to the semisimplicity of the $A_j$, it is also the Lie algebra of the unipotent radical of the tannakian group of the object $M_1\oplus M_2$ of $\bT$. Since $\fu(M_1)$ and $\fu(M_2)$ are semisimple objects of $\bT$, the restriction map 
\[
\bigoplus_{j=1}^2 \Hom(\fu(M_j), \inHom(A_3,A_1)) \ \rightarrow \ \Hom(\fu(M_1\oplus M_2), \inHom(A_3,A_1))
\]
is surjective. Moreover, since $\fu(M_1)$ and $\fu(M_2)$ are abelian, so is $\fu(M_1\oplus M_2)$. In view of Proposition \ref{prop: map for Ext groups}, we thus obtain a surjection
\[
\bigoplus_{j=1}^2 \Ext^1_{\langle M_j\rangle^{\otimes}}(A_3,A_1) \twoheadrightarrow 
\Ext^1_{\langle M_1\oplus M_2\rangle^{\otimes}}(A_3,A_1).
\]
Since $\fu_{-2}(M_1)$ and $\fu_{-2}(M_2)$ are zero, by Corollary \ref{cor: when is [u,u]=u-2} the two Ext groups on the left vanish. We thus get
\[
\Ext^1_{\langle M_1\oplus M_2\rangle^{\otimes}}(A_3,A_1) = 0.
\]
This forces $M_1$ and $M_2$ to be the same in $\Extpan(N,L)$. Indeed, the isomorphism classes of $M_1$ and $M_2$ (as blended extensions) belong to the subset $\Extpan_{\langle M_1\oplus M_2\rangle^{\otimes}}(N,L)$ of $\Extpan(N,L)$ consisting of the isomorphism classes of blended extensions of $N$ by $L$ in the category $\langle M_1\oplus M_2\rangle^{\otimes}$. The set $\Extpan_{\langle M_1\oplus M_2\rangle^{\otimes}}(N,L)$ is a torsor for $\Ext^1_{\langle M_1\oplus M_2\rangle^{\otimes}}(A_3,A_1)$. Since this Ext group vanishes, $\Extpan_{\langle M_1\oplus M_2\rangle^{\otimes}}(N,L)$ is a singleton.
\end{proof}

By Proposition \ref{prop: [u,u]={u-1,u-1}}, an obvious necessary condition for existence of a blended extension of $N$ by $L$ with vanishing $\fu_{-2}$ is that $\{\fu_{-1},\fu_{-1}\}$ must vanish (where $\{,\}$ is the pairing of \eqref{eq21}). When $\bT$ is the category of differential modules over a differential field of characteristic zero with an algebraically closed constant field, a result of Hardouin \cite[Th\'{e}or\`{e}me 1.1(i)]{Har05} asserts that the condition $\{\fu_{-1},\fu_{-1}\}=0$ is also sufficient for existence of a blended extension with vanishing $\fu_{-2}$ (assuming $\Extpan(N,L)$ is nonempty). The same should be true in general, but we have not tried to prove it. For arbitrary $\bT$, Bertrand has proved in \cite{Ber13} that in a rather special setting where he defines a notion of self-duality for blended extensions, assuming $\Extpan(N,L)$ is nonempty, there always exists a blended extension $M$ of $N$ by $L$ such that $\fu_{-2}(M)=[\fu(M),\fu(M)]$. (See {\it loc. cit.}, Theorems 1 and 2, also the proof of the latter.) The same might also be true in general.

\section{Application to the Hodge conjecture for 1-motives}\label{sec: application to HN for 1-motives}

Throughout this section $\FF$ is an algebraically closed subfield of $\CC$.
\subsection{} Let $\mathbf{NMM}(\mathbb{F})$ be Nori's tannakian category of (ineffective mixed) motives over $\FF$ \cite{HM17}. We refer to the objects of $\mathbf{NMM}(\mathbb{F})$ as Nori motives. Let $\mathbf{MHS}$ be the category of rational mixed Hodge structures. In \cite{An19} Andr\'{e} proves the following theorem:
\begin{thm}[Andr\'{e}]\label{thm: Andre's thm}
Let $M$ be a 1-motive\footnote{Note the slight change in notation in this section; we will be using $M$ for a Deligne 1-motive, not the motive associated with it in a tannakian category of motives.} over $\FF$ (in the sense of Deligne \cite{De74}). Then the motivic Galois group of the Nori motive associated with $M$ coincides with its Mumford-Tate group.
\end{thm}
Here, by the motivic Galois group of a Nori motive one means the fundamental group of the Nori motive with respect to the ``Betti" fiber functor, which is the composition of the Hodge realization functor $\mathbf{NMM}(\mathbb{F}) \rightarrow \mathbf{MHS}$ and the forgetful functor from $\mathbf{MHS}$ to the category of finite-dimensional rational vector spaces. The Mumford-Tate group of the motive is the fundamental group of its Hodge realization with respect to the forgetful fiber functor. The Mumford-Tate group is always canonically a subgroup of the motivic Galois group, and the motivic version of the Hodge conjecture for Nori motives predicts that the two groups are the same for every object of $\mathbf{NMM}(\mathbb{F})$. 

Andr\'{e}'s proof of Theorem \ref{thm: Andre's thm} uses a deformation argument to reduce the problem to the case of semisimple 1-motives. This semisimple case was proved earlier by Andr\'{e} himself in \cite{An96} in the setting of motives via motivated correspondences. Arapura \cite{Ar13} has proved that Andr\'{e}'s category of motives constructed via motivated correspondences is canonically equivalent to the full subcategory of $\mathbf{NMM}(\mathbb{F})$ consisting of semisimple objects, so that Andr\'{e}'s result in \cite{An96} also resolves the case of semisimple 1-motives for the Nori setting.

In \cite{An19} Andr\'{e} asks if it is possible to give a proof of the reduction to the semisimple case that is not based on a deformation argument. Here, as an application of the earlier results of the paper, we propose such an argument. The argument is rather formal, and works in any tannakian category of motives as long as the Hodge realization map is injective on extensions of $\mathbbm{1}$ by $\QQ(1)$, as we explain below.

\subsection{}
Let $\mathbf{DM}_1(\FF)$ be Deligne's abelian category of 1-motives over $\FF$ up to isogeny. Suppose that $\mathbf{MM}(\FF)$ is a tannakian category over $\QQ$ of mixed motives over $\FF$, with a $\QQ$-linear (covariant) functor
\[
h: \mathbf{DM}_1(\FF) \rightarrow \mathbf{MM}(\FF)
\]
the composition of which with the Hodge realization functor $\fR^{H}: \mathbf{MM}(\FF)\rightarrow \mathbf{MHS}$ (which is exact, faithful, linear, respectful of the weight filtration and tensor products) is the Hodge realization functor on 1-motives constructed by Deligne in \cite{De74} (thus in particular, $h$ is faithful and exact). Suppose moreover that pure objects of $\mathbf{MM}(\FF)$ are semisimple, and that $h$ behaves well with respect to duals, i.e., for every 1-motive $M$ we have a canonical isomorphism $h(M^\ast)\cong h(M)^\vee(1)$ functorial in $M$, where $M^\ast$ is the Cartier dual of $M$ and $-^\vee$ is the dualizing functor in $\mathbf{MM}(\FF)$ in the tannakian sense. These requirements are rather basic, and any of the known constructions of a tannakian category of mixed motives (in particular, Nori's $\mathbf{NMM}(\FF)$) satisfies these conditions.

Given a 1-motive $M$ over $\FF$, thanks to its weight filtration, $M$ fits into a blended extension as in \eqref{eq1} in the abelian category $\mathbf{DM}_1(\FF)$. The top left object, the top right object, and the bottom right object (in places of $A_1$, $A_2$, and $A_3$) respectively are the torus part, the abelian part, and the lattice part of $M$. Applying $h$ and then further $\fR^H$ to this blended extension we obtain blended extensions $h(M)$ and $\fR^H h(M)$ in the tannakian categories $\mathbf{MM}(\FF)$ and $\mathbf{MHS}$. Following the constructions of \S \ref{sec: initial considerations}, we thus may speak of $\fu(h(M))$, $\fu_{-1}(h(M))$ and $\fu_{-2}(h(M))$, as well as the corresponding objects $\fu(\fR^{H} h(M))$, $\fu_{-1}(\fR^{H} h(M))$ and $\fu_{-2}(\fR^{H} h(M))$ for the Hodge realization of $h(M)$. The latter three objects, respectively, are canonically contained in the Hodge realizations of the former three objects. In fact, the top extension in diagram \eqref{eq10} of \S \ref{sec: def of u-1 and u-2} for $\fR^H h(M)$ is contained in the Hodge realization of the analogous extension for $h(M)$.

Since $\FF$ is algebraically closed and the composition
\[
\mathbf{DM}_1(\FF) \xrightarrow{h} \mathbf{MM}(\FF) \xrightarrow{\fR^H} \mathbf{MHS}
\]
coincides with the usual Hodge realization of 1-motives, the composition is a full functor (see \cite[Proposition 2.1]{An19}). Since $\fR^H$ is faithful, it follows that the functor $h$ is also full. We thus have an injection
\[
\Ext^1_{\mathbf{DM}_1(\FF)}(M,M') \hookrightarrow \Ext^1_{\mathbf{MM}(\FF)}(h(M),h(M'))
\]
for every 1-motives $M$ and $M'$ over $\FF$. Deligne conjectures that in a good tannakian category of mixed motives the essential image of $\mathbf{DM}_1(\FF)$ should be closed under extensions (\cite{De89}, \S 2.4). Thus the above map should be an isomorphism. In particular, taking $M=\ZZ$ and $M'=\mathbb{G}_m$, we should have that
\begin{equation}\label{eq26}
\Ext^1_{\mathbf{MM}(\FF)}(\mathbbm{1},\QQ(1)) \cong \Ext^1_{\mathbf{DM}_1(\FF)}(\ZZ,\mathbb{G}_m).
\end{equation}
By the fullness of $\fR^H h$, the composition
\[
\Ext^1_{\mathbf{DM}_1(\FF)}(\ZZ,\mathbb{G}_m) \xrightarrow{h} \Ext^1_{\mathbf{MM}(\FF)}(\mathbbm{1},\QQ(1)) \xrightarrow{ \ \fR^{H} \ } \Ext^1_{\mathbf{MHS}}(\mathbbm{1},\QQ(1))
\]
is injective. (In fact, identifying $\Ext^1_{\mathbf{DM}_1(\FF)}(\ZZ,\mathbb{G}_m)$ with $\FF^\times\otimes \QQ$ and $\Ext^1_{\mathbf{MHS}}(\mathbbm{1},\QQ(1))$ with $\CC/2\pi i\QQ$, the composition above is just given by the logarithm function.) Thus the special case \eqref{eq26} of Deligne's conjecture would imply that the Hodge realization map on $\Ext^1_{\mathbf{MM}(\FF)}(\mathbbm{1},\QQ(1))$ is injective.

We shall prove the following:

\begin{thm}\label{thm: application to Hodge-Nori}
Assume that $\mathbf{MM}(\FF)$ is a tannakian category of mixed motives\footnote{To be clear, in addition to the hypotheses explicitly mentioned in the statement, all that here is needed from $\mathbf{MM}(\FF)$ is the following: $\mathbf{MM}(\FF)$ is a tannakian category over $\QQ$ with a weight filtration similar to the weight filtration on mixed Hodge structures, there is an exact faithful tensor $\QQ$-linear functor $\fR^H: \mathbf{MM}(\FF)\rightarrow \mathbf{MHS}$ preserving the weight filtration, and the images of abelian varieties under $h$ are semisimple. Working in this generality, the object $\QQ(1)$ of $\mathbf{MM}(\FF)$ is defined to be $h(\mathbb{G}_m)$.} over $\FF$ with a $\QQ$-linear functor $h: \mathbf{DM}_1(\FF) \rightarrow \mathbf{MM}(\FF)$ which behaves well with respect to duals (in the sense mentioned above), and whose composition with the Hodge realization $\fR^{H}: \mathbf{MM}(\FF)\rightarrow \mathbf{MHS}$ is the usual Hodge realization of 1-motives constructed by Deligne. Let $M$ be a 1-motive over $\FF$.

\begin{itemize}[wide]
\item[(a)] Then
\[
\fu_{-1}(\fR^{H} h(M)) = \fR^{H} \fu_{-1}(h(M)).
\]
\item[(b)] Suppose that the Hodge realization map 
\[
\Ext^1_{\mathbf{MM}(\FF)}(\mathbbm{1},\QQ(1))\rightarrow \Ext^1_{\mathbf{MHS}}(\mathbbm{1},\QQ(1))
\]
is injective on the subgroup $\Ext^1_{\langle h(M)\rangle^{\otimes}}(\mathbbm{1},\QQ(1))$ of $\Ext^1_{\mathbf{MM}(\FF)}(\mathbbm{1},\QQ(1))$. Then 
\begin{equation}\label{eq44}
\fu(\fR^{H} h(M)) = \fR^{H} \fu(h(M)).
\end{equation}
\end{itemize}
\end{thm}

\begin{proof}
Let $M=[\ZZ^n\xrightarrow{v}G]$, where $G$ is a semiabelian variety over $\FF$, an extension of an abelian variety $A$ by a torus $\mathbb{G}_m^s$. Consider the blended extension in $\mathbf{DM}_1(\FF)$ given by the weight filtration on $M$. Its top row is the extension $G$ of $A$ by $\mathbb{G}_m^s$, and its right column is the extension $M/\mathbb{G}_m^s=[\ZZ^n\rightarrow A]$ (the map being the projection of $v$) of $\ZZ^n$ by $A$. The blended extensions to which the earlier results of the paper will be applied are $h(M)$ and $\fR^Hh(M)$. Both of these blended extensions do indeed satisfy conditions (C1) and (C2) of \S \ref{sec: conditions C1 and C2}. The equality \eqref{eq44} holds if and only if we have equalities at the levels of both $\fu_{-1}$'s and $\fu_{-2}$'s.

\begin{itemize}[wide]
\item[(a)] We will show the equality of $\fR^{H}\fu_{-1}(h(M))$ and $\fu_{-1}(\fR^{H} h(M))$ by comparing their descriptions given by Theorem \ref{thm: characterization of u-1}. Set $N=h(M/\mathbb{G}_m^s)$ and $L=h(G)$, the former an extension of $\mathbbm{1}^n$ by $h(A)$ and the latter an extension of $h(A)$ by $\QQ(1)^s$. Consistent with the notation used in \S \ref{sec: char of u-1}, the corresponding extensions of $\mathbbm{1}$ by 
\[
\inHom(\mathbbm{1}^n, h(A)) \cong h(A^n) \hspace{.2in}\text{and}\hspace{.2in} \inHom(h(A),\QQ(1)^s)\stackrel{(\dagger)}{\cong}  h({A^\ast}^s)
\]
will be respectively denoted by $\sN$ and $\sL$. The identification $(\dagger)$ here uses the canonical isomorphism $h(A^\ast)\cong h(A)^\vee(1)$. Using the functoriality of the isomorphisms $h(-^\ast)\cong h(-)^\vee(1)$, it is not difficult to see that up to a sign, the extension $\sL$ in $\Ext^1_{\mathbf{MM}(\FF)}(\mathbbm{1}, h({A^\ast}^s))$ coincides with the element associated with $h(G^\ast)$.

We have an element
\[
(\sL,\sN) \in \Ext^1_{\mathbf{MM}(\FF)}(\mathbbm{1}, h({A^\ast}^s)\oplus h(A^n) ) \cong \Ext^1_{\mathbf{MM}(\FF)}(\mathbbm{1}, h({A^\ast}^s\times A^n)).  
\]
By Theorem \ref{thm: characterization of u-1}, the subobject $\fu_{-1}(h(M))$ of $h({A^\ast}^s\times A^n)$ is the intersection of the kernels of all the endomorphisms of $h({A^\ast}^s\times A^n)$  which annihilate $(\sL,\sN)$. Applying $\fR^H$ to this, we get
\[\fR^H \fu_{-1}(h(M)) \ \ \ = \bigcap_{\substack{\phi\in \End(h({A^\ast}^s\times A^n)) \\ \phi_\ast(\sL,\sN)=0}} \ker(\fR^H \phi).\]
On the other hand, applying Theorem \ref{thm: characterization of u-1} to $\fR^H h(M)$ we have
\[\fu_{-1}(\fR^H h(M)) \ \ \ = \bigcap_{\substack{\psi\in \End(\fR^H h({A^\ast}^s\times A^n)) \\ \psi_\ast(\fR^H \sL,\fR^H\sN)=0}} \ker(\psi).\]
It thus suffices to show that every endomorphism $\psi$ of $\fR^H h({A^\ast}^s\times A^n)$ which annihilates  $(\fR^H \sL,\fR^H\sN)$ is of the form $\fR^H \phi$ for some (unique) endomorphism $\phi$ of $h({A^\ast}^s\times A^n)$ that annihilates $(\sL,\sN)$. 

By fullness of $\fR^H h$, we have isomorphisms
\begin{equation}\label{eq32}
\End_{\mathbf{DM}_1(\FF)}({A^\ast}^s\times A^n) \xrightarrow{ \ h \ } \End(h({A^\ast}^s\times A^n)) \xrightarrow{ \, \fR^H \, } \End(\fR^H h({A^\ast}^s\times A^n)).
\end{equation}
Thus it suffices to argue that for every $\alpha\in \End_{\mathbf{DM}_1(\FF)}({A^\ast}^s\times A^n)$, we have an implication
\begin{equation}\label{eq33}
(\fR^H h(\alpha))_\ast(\fR^H\sL,\fR^H\sN)=0 \ \ \ \Rightarrow \ \ \ (h(\alpha))_\ast(\sL,\sN)=0.
\end{equation}
There is a commutative diagram 
\begin{equation}\label{eq28} 
\begin{tikzcd}[row sep = large]
\Ext^1_{\mathbf{DM}_1(\FF)}(\ZZ, {A^\ast}^s\times A^n) \ar[r, hookrightarrow, "h"]  \ar[dr, hookrightarrow, "\fR^H h " ] & \Ext^1_{\mathbf{MM}(\FF)}(\mathbbm{1}, h({A^\ast}^s\times A^n)) \ar[d, "\fR^H"] \\
& \Ext^1_{\mathbf{MHS}}(\mathbbm{1}, \fR^H h({A^\ast}^s\times A^n)).
\end{tikzcd}
\end{equation}
Each endomorphism algebra in \eqref{eq32} acts by pushforwards on the corresponding $\Ext^1$ group above. Moreover, the maps in \eqref{eq28} commute with the actions of these endomorphism algebras (as $h$ and $\fR^H$ are exact). The implication \eqref{eq33} now follows, since the extension $(\sL,\sN)$ is in the image of $\Ext^1_{\mathbf{DM}_1(\FF)}(\ZZ, {A^\ast}^s\times A^n)$.

\medskip\par 
\item[(b)] We need to show that
\begin{equation}\label{eq34}
\fu_{-2}(\fR^{H} h(M)) = \fR^{H} \fu_{-2}(h(M)).
\end{equation}
By the previous part and Proposition \ref{prop: [u,u]={u-1,u-1}} (and the fact that $\fR^H$ takes the bracket $\{  ,  \}$ of $h(M)$ to the counterpart for $\fR^H h(M)$), we have
\[
[\fu(\fR^H h(M)),\fu(\fR^H h(M))] = \fR^H[\fu(h(M)),\fu(h(M))].
\]
Thus we have an injection
\[
\frac{\fu_{-2}(\fR^H h(M))}{[\fu(\fR^H h(M)),\fu(\fR^H h(M))]} \stackrel{i}{\hookrightarrow} \frac{\fR^H \fu_{-2}(h(M))}{\fR^H[\fu(h(M)),\fu(h(M))]}.
\]
The isomorphisms of Proposition \ref{prop: characterization of u-2/[u,u]} (or rather, its proof with $\QQ(1)$ replacing $\inHom(A_3,A_1)$) for $h(M)$ and $\fR^H h(M)$ give rise to a commutative diagram
\[
\begin{tikzcd}
\Ext^1_{\langle h(M)\rangle^{\otimes}}(\mathbbm{1},\QQ(1)) \arrow[dd, " \, \fR^H" '] \arrow[r, "\simeq "] & \Hom\bigm(\displaystyle{\frac{\fu_{-2}(h(M))}{[\fu (h(M)),\fu (h(M))]}}, \QQ(1)\bigm) \arrow[d, "\fR^H \, " ' , "\simeq" ]\\
& \Hom\bigm(\displaystyle{\frac{\fR^H\fu_{-2}( h(M))}{\fR^H [\fu (h(M)),\fu (h(M))]}}, \QQ(1)\bigm)\arrow[d, twoheadrightarrow, "i^\ast \, " ' ] \\
\Ext^1_{\langle \fR^H h(M)\rangle^{\otimes}}(\mathbbm{1},\QQ(1)) \arrow[r, "\simeq"] & \Hom\bigm(\displaystyle{\frac{\fu_{-2}(\fR^H h(M))}{[\fu(\fR^H h(M)),\fu(\fR^H h(M))]}}, \QQ(1)\bigm).
\end{tikzcd}
\]
The Hodge realization map on the right is an isomorphism as $\fu_{-2}(h(M))/[\fu(h(M)),\fu(h(M))]$ is a subquotient of $\inHom(\mathbbm{1}^n,\QQ(1)^s)$ and hence is a direct sum of copies of $\QQ(1)$. The pullback map $i^\ast$ is surjective because $i$ admits a retraction (by semisimplicity).  The Hodge realization map on the left is crucially assumed to be injective in the statement of the theorem. It follows that $i^\ast$ is injective, so that
\[
\Hom(\frac{\fR^H \fu_{-2}(h(M))}{\fu_{-2}(\fR^H h(M))},\QQ(1)) = 0.
\]
But $\displaystyle{\frac{\fR^H \fu_{-2}(h(M))}{\fu_{-2}(\fR^H h(M))}}$ is a direct sum of copies of $\QQ(1)$, so that we get \eqref{eq34}.
\end{itemize}
\end{proof}

\begin{rem}
Prior to \cite{An19}, Jossen had proved in \cite[Theorem 6.2]{Jo14} that for 1-motives over $\CC$, the unipotent radicals of the motivic Galois and Mumford-Tate groups coincide. (See  also the appendix of the same paper.)
\end{rem}


\begin{thebibliography}{00}


\bibitem{An96}
Y. Andr\'{e}, 
Pour une th\'{e}orie inconditionelle des motifs, 
Publ. Math. I. H. E. S. 83 (1996) 5-49

\bibitem{An04}
Y. Andr\'{e},
Une introduction aux motifs,
Soci\'{e}t\'{e} Math\'{e}matique de France, 2004


\bibitem{An19} 
Y. Andr\'{e},
A note on 1-motives,
International Mathematics Research Notices, 2021, no. 3, 2074-2080

\bibitem{Ar13}
D. Arapura, 
An abelian category of motivic sheaves, 
Adv. Math. 233 1 (2013), 135-195


\bibitem{ABV15}
J. Ayoub and L. Barbieri-Viale,
Nori 1-motives,
Math. Ann., Vol. 361 (2015), 367-402


\bibitem{Be02}
C. Bertolin,
The Mumford-Tate group of 1-motives,
Ann. Inst. Fourier, Grenoble, 52, 4 (2002), 1041-1059

\bibitem{Be03}
C. Bertolin,
Le radical unipotent du groupe de Galois motovique d'un 1-motif,
Math. Ann. 327, 585-607 (2003)


\bibitem{BPSS}
C. Bertolin, P. Philippon, B. Saha, E. Saha,
Semi-abelian analogues of Schanuel conjecture and applications,
Journal of Algebra,
Volume 596,
2022,
Pages 250-288

\bibitem{BP1}
C. Bertolin and P. Philippon,
Mumford-Tate groups of 1-motives and Weil pairing,
J. Pure Appl. Algebra, Volume 228, Issue 10, 2024, 107702

\bibitem{Ber01}
D. Bertrand,
Unipotent radicals of differential Galois group and integrals of solutions of inhomogeneous equations,
Math. Ann. 321 (2001), no. 3, 645–666

\bibitem{Ber13}
D. Bertrand,
Extensions panach\'{e}es autoduales,
J. K-Theory 11 (2013), no. 2, 393–411

\bibitem{De74}
P. Deligne,
Theorie de Hodge III,
Publications Math\'{e}matiques de l'I.H.\'{E}.S., tome 44 (1974), p. 5-77

\bibitem{DM82}
P. Deligne and J.S. Milne,
Tannakian Categories,
In Hodge Cycles, Motives, and Shimura Varieties,
Lecture Notes in Mathematics 900, Springer-Verlog, Berlin (1982)

\bibitem{De89}
P. Deligne,
Le groupe fondamental de la droite projective moins trois points,
Galois Groups over $\QQ$, MSRI Publ. 16, pp. 79-313, Springer-Verlag (1989)

\bibitem{De90}
P. Deligne,
Cat\'{e}gories tannakiennes,
Grothendieck Festschrift Vol II, Progress in Mathematics, 87 ( Birkhäuser Boston 1990) pp. 111–195

\bibitem{EM1}
P. Eskandari and V. K. Murty,
The fundamental group of an extension in a Tannakian category and the unipotent radical of the Mumford-Tate group of an open curve, Pacific Journal of Mathematics, Vol. 325 (2023), No. 2, 255-279

\bibitem{EM2}
P. Eskandari and V. K. Murty,
On unipotent radicals of motivic Galois groups,
Algebra \& Number Theory 17-1 (2023), pp 165-215

\bibitem{Es23}
P. Eskandari,
On blended extensions in filtered abelian categories and motives with maximal unipotent radicals, Documenta Mathematica (2025), published online first, DOI 10.4171/DM/1052 

\bibitem{Es26}
P. Eskandari,
Depth one part of Tannakian groups of filtrations, 
Canadian Journal of Mathematics (2026), pp. 1-30. doi:10.4153/S0008414X26102168.

\bibitem{Gr68}
A. Grothendieck,
Mod\`{e}les de N\'{e}ron et monodromie,
SGA VII.1, no 9, Springer LN 288, 1968


\bibitem{Har05}
C. Hardouin, 
Calcul du groupe de Galois du produit de  trois op\'{e}rateurs diff\'{e}rentiels compl\`{e}tement r\'{e}ductibles
C. R. Acad. Sci. Paris, Ser. I 341 (2005)


\bibitem{Har06}
C. Hardouin,
Hypertranscendance et Groupes de Galois aux diff\'{e}rences, arXiv 0609646v2, 2006


\bibitem{Har11}
C. Hardouin,
Unipotent radicals of Tannakian Galois groups in positive characteristic, 
Arithmetic and Galois theories of differential equations, 223-239, S\'{e}min. Congr., 23, Soc. Math. France, Paris, 2011

\bibitem{HM17}
A. Huber and S. M\"{u}ller-Stach,
Periods and Nori Motives,
Springer, 2017


\bibitem{Jo14}
P. Jossen,
On the Mumford-Tate conjecture for 1-motives,
Inventiones Math. (2014) 195: 393-439


\end{thebibliography}
\end{document}